\pdfoutput=1
\RequirePackage{silence}
\documentclass[a4paper,svgnames]{amsart}
\usepackage[hmarginratio={1:1},vmarginratio={1:1},lmargin=60.0pt,tmargin=60.0pt]{geometry}
\usepackage{float}

\usepackage{booktabs}
\allowdisplaybreaks

\usepackage[T1]{fontenc}
\usepackage{latexsym,exscale,amsfonts,amssymb,mathtools}
\usepackage{amsmath,amsthm,amsfonts,amssymb,amscd,textcomp,bbm}
\usepackage{mathrsfs,stackrel}
\usepackage{etoolbox}

\usepackage[all]{xy}
\usepackage{tikz}
\usetikzlibrary{cd}
\usetikzlibrary{decorations}
\usetikzlibrary{decorations.markings}
\usetikzlibrary{decorations.pathreplacing}
\usetikzlibrary{decorations.pathmorphing}
\usetikzlibrary{arrows.meta,shapes,positioning,matrix,calc}
\usetikzlibrary{shapes.callouts}
\tikzset{anchorbase/.style={baseline={([yshift=-0.5ex]current bounding box.center)}},
tinynodes/.style={font=\tiny, text height=0.25ex, text depth=0.05ex},
smallnodes/.style={font=\scriptsize, text height=0.75ex, text depth=0.15ex},
crossline/.style={preaction={draw=white,line width=5.0pt,-},preaction={draw=black,line width=0.9pt,-}},
usual/.style={line width=1.0,color=black},
dot/.style = {
decoration={markings,
post length=0.25mm,
pre length=0.25mm,
mark=at position #1 with {\node[circle,radius=0.2cm,inner sep=-1.5pt,color=black,fill=black]{};}
},
postaction={decorate}
},
dot/.default=1,
}
\usepackage{tikz}
\usetikzlibrary{cd}
\usetikzlibrary{decorations}
\usetikzlibrary{decorations.markings}
\usetikzlibrary{decorations.pathreplacing}
\usetikzlibrary{decorations.pathmorphing}
\usetikzlibrary{arrows.meta,shapes,positioning,matrix,calc}
\usetikzlibrary{shapes.callouts}
\tikzset{anchorbase/.style={baseline={([yshift=-0.5ex]current bounding box.center)}},
tinynodes/.style={font=\tiny, text height=0.25ex, text depth=0.05ex},
smallnodes/.style={font=\scriptsize, text height=0.75ex, text depth=0.15ex},
crossline/.style={preaction={draw=white,line width=5.0pt,-},preaction={draw=black,line width=0.9pt,-}},
usual/.style={line width=1.0,color=black},
dot/.style = {
decoration={markings,
post length=0.25mm,
pre length=0.25mm,
mark=at position #1 with {\node[circle,radius=0.2cm,inner sep=-1.5pt,color=black,fill=black]{};}
},
postaction={decorate}
},
dot/.default=1,
}
\usepackage{aliascnt}
\usepackage{ytableau}
\usepackage{xparse}
\usepackage{cite}

\usepackage{dynkin-diagrams}
\usepackage{nicematrix}
\usepackage{tikz-3dplot}
\usepackage{mathrsfs}
\DeclareMathAlphabet{\mathscrbf}{OMS}{mdugm}{b}{n}
\DeclareMathOperator{\M}{M}
\DeclareMathOperator{\Ext}{Ext}
\DeclareMathOperator{\GL}{GL}

\usepackage{array}
\newcolumntype{C}{>{$}c<{$}}
\newcolumntype{P}[1]{>{\centering\arraybackslash}p{#1}} 
\DeclarePairedDelimiter\floor{\lfloor}{\rfloor}

\usepackage{xcolor,colortbl}

\definecolor{mygray}{gray}{0.6}
\definecolor{mygraydark}{gray}{0.4}
\definecolor{mygraylight}{gray}{0.85}
\definecolor{spinach}{RGB}{46,139,87}
\definecolor{tomato}{RGB}{255,99,71}
\definecolor{orchid}{RGB}{143,40,194}
\definecolor{neon}{RGB}{77,77,255}
\definecolor{lightneon}{RGB}{110,110,255}
\definecolor{pumpkin}{RGB}{224,180,80}
\definecolor{citron}{RGB}{190,180,90}

\definecolor{lava}{RGB}{207,16,32}
\definecolor{cream}{RGB}{255,253,208}
\definecolor{verdigris}{RGB}{67,179,174}
\definecolor{Black}{RGB}{0,0,0}
\definecolor{mydarkblue}{RGB}{10,10,170}
\definecolor{darkspinach}{RGB}{20,70,20}
\definecolor{darktomato}{RGB}{155,40,30}
\definecolor{darkorchid}{RGB}{50,10,100}
\definecolor{darklava}{RGB}{150,8,16}

\definecolor{zero}{RGB}{0,0,0}
\definecolor{one}{RGB}{255,0,0}
\definecolor{two}{RGB}{0,255,0}
\definecolor{three}{RGB}{0,0,255}

\usepackage{enumitem}
\setlist[enumerate]{itemsep=0.15cm,label=\emph{\upshape(\alph*)}}
\setlist[enumerate,2]{itemsep=0.15cm,label=\emph{\upshape(\roman*)}}
\setlist[enumerate,3]{itemsep=0.15cm,label=\emph{\upshape(\Alph*)}}

\newcommand{\N}{\mathbb{N}}
\newcommand{\Z}{\mathbb{Z}}

\newcommand{\C}{\mathbb{C}}

\DeclareMathOperator{\Ind}{Ind}
\DeclareMathOperator{\Res}{Res}

\DeclareMathOperator{\SLM}{SLM}

\DeclareMathOperator{\SL}{SL}
\newcommand{\ceil}[1]{\left\lceil #1 \right\rceil}

\def\NewTheorem#1{%
\newaliascnt{#1}{equation}%
\newtheorem{#1}[#1]{#1}%
\aliascntresetthe{#1}%
\expandafter\def\csname #1autorefname\endcsname{#1}%
}
\def\equationautorefname~#1\null{(#1)\null}

\numberwithin{equation}{subsection}

\NewTheorem{Proposition}
\NewTheorem{Theorem}
\NewTheorem{Corollary}
\AtEndEnvironment{Corollary}{\null\hfill$\square$}%
\NewTheorem{Lemma}
\NewTheorem{Conjecture}
\NewTheorem{Speculation}
\NewTheorem{Observation}
\NewTheorem{Assumption}
\theoremstyle{definition}
\NewTheorem{Definition}
\AtEndEnvironment{Definition}{\null\hfill$\Diamond$}%
\NewTheorem{Classification Problem}
\AtEndEnvironment{Classification Problem}{\null\hfill$\Diamond$}%
\NewTheorem{Notation}
\AtEndEnvironment{Notation}{\null\hfill$\Diamond$}%
\NewTheorem{Example}
\AtEndEnvironment{Example}{\null\hfill$\Diamond$}%
\NewTheorem{Examples}
\AtEndEnvironment{Examples}{\vskip-10mm\null\hfill$\Diamond$}%

\theoremstyle{remark}
\NewTheorem{Remark}
\AtEndEnvironment{Remark}{\null\hfill$\Diamond$}%
\NewTheorem{Question}
\AtEndEnvironment{Question}{\null\hfill$\Diamond$}%

\usepackage[hypertexnames=false]{hyperref}
\usepackage{bookmark}
\hypersetup{
pdftoolbar=true,
pdfmenubar=true,
pdffitwindow=false,
pdfstartview={FitH},
pdftitle={Tensor powers of representations of (diagram) monoids},
pdfauthor={David He and Daniel Tubbenhauer},
pdfsubject={},
pdfcreator={David He and Daniel Tubbenhauer},
pdfproducer={David He and Daniel Tubbenhauer},
pdfkeywords={},
pdfnewwindow=true,
colorlinks=true,
linkcolor=mydarkblue,
citecolor=teal,
filecolor=magenta,
urlcolor=orchid,
linkbordercolor=lava,
citebordercolor=teal,
urlbordercolor=orchid,
linktocpage=true
}

\def\makeautorefname#1#2{\csdef{#1autorefname}{#2}}
\makeautorefname{section}{Section}%
\makeautorefname{subsection}{Section}%
\makeautorefname{subsubsection}{Section}%

\begin{document}

\arrayrulewidth=0.5mm
\setlength{\arrayrulewidth}{0.5mm}

\title{Tensor powers of representations of (diagram) monoids}
\author{David He and Daniel Tubbenhauer}
\address{D.H.: The University of Sydney, School of Mathematics and Statistics, Australia}
\email{d.he@sydney.edu.au}

\address{D.T.: The University of Sydney, School of Mathematics and Statistics F07, Office Carslaw 827, NSW 2006, Australia, \href{http://www.dtubbenhauer.com}{www.dtubbenhauer.com}, \href{https://orcid.org/0000-0001-7265-5047}{ORCID 0000-0001-7265-5047}}
\email{daniel.tubbenhauer@sydney.edu.au}

\begin{abstract}
We study tensor powers of representations of finite monoids, focusing on the growth behavior of their composition length and the number of indecomposable summands. Special attention is given to diagram monoids such as the Temperley--Lieb, Motzkin, and Brauer monoids. For these examples, we compute explicit data, including some character tables, and analyze patterns in the decomposition of their tensor powers.
\end{abstract}


\addtocontents{toc}{\protect\setcounter{tocdepth}{1}}

\maketitle

\tableofcontents

\section{Introduction}

Throughout let $M$ be a finite monoid, with group of units $G$, and let $k$ be an algebraically closed field of characteristic $p\ge 0$. We associate to a $kM$-module $V$ two statistics, denoted by $l(n)$ and $b(n)$:
\[l(n)=l^{M,V}(n)= \text{length of\ } V^{\otimes n}. \]
\[b(n)=b^{M,V}(n)=\# M\text{-indecomposable summands of\ } V^{\otimes n} \text{(counted with multiplicity)}.\]

\begin{Remark}
In general, the study of $b(n)$ is more difficult than that of $l(n)$ (while $b(n)$ counts summands in an honest decomposition of $V^{\otimes n}$, $l(n)$ counts summands in the decomposition of the Grothendieck class of $V^{\otimes n}$ and $l(n)$ can be seen as the semisimplification of $b(n)$). For example, a general exact formula for $b(n)$ appears impossible already in the case when $M$ is a group. Of course, $b(n)=l(n)$ when $kM$ is semisimple, but in general we only have $b(n)\leq l(n)$. 
\end{Remark}

The growth problem for $V$ is concerned with  understanding (the asymptotics of) the sequences $l(n)$ and $b(n)$. Concretely, it involves the following set of questions:
\begin{enumerate}[label={(\arabic*)}]
\item \label{1} Can we find a formula for $l(n)$ (resp. $b(n)$), or more generally, a formula for an asymptotic expression $k(n)$ with $l(n)\sim k(n)$ (resp. $a(n)$ such that $b(n)\sim a(n)$)? Here $l(n)\sim k(n)$ if they are asymptotically equal: ${l(n)}/{k(n)}\xrightarrow{ n \to \infty } 1$, and ditto for $a(n)$ and $b(n)$.
\item \label{2} Can we understand the \textit{rate of convergence} by quantifying how fast $\lvert l(n)/k(n)-1\rvert$ (resp. $\lvert b(n)/a(n)-1\rvert$) converges to 0?
\item \label{3} Similarly, can we bound the \textit{variance} $\lvert l(n)-k(n)\rvert$ (resp. $\lvert b(n)-a(n)\rvert$)?
\end{enumerate}
If $S$ is a simple (resp. indecomposable) $kM$-module, we are sometimes also interested in the finer statistics: 
\begin{enumerate}[label={(\arabic*)},resume]
\item What is the composition factor multiplicity $[V^{\otimes n}:S]$ of $S$ in $V^{\otimes n}$ (resp. $[V^{\otimes n}:S]_b$, the number of copies of $S$ in the Krull--Schmidt decomposition of $V^{\otimes n}$) and its asymptotics?
\end{enumerate}

\begin{Remark}
More generally, growth problems may be defined for any additive Krull--Schmidt monoidal category (see \cite{lacabanne2024asymptotics}). The growth problems and related questions for various categories have been studied in many papers, for recent examples we refer the reader to \cite{coulembier2023asymptotic, coulembier2024fractalbehaviortensorpowers,Coulembier_2023,he2024growthproblemsrepresentationsfinite, he2025growthproblemsrepresentationsfinite,lacabanne2023asymptotics, lacabanne2024asymptotics,lachowska2024tensorpowersvectorrepresentation, larsen2024boundsmathrmsl2indecomposablestensorpowers, gruber2025growthproblemsdiagramcategories}. 
\end{Remark}

In \cite[Theorem 3B.2]{he2025growthproblemsrepresentationsfinite}, the $l(n)$ version of the growth problem is answered in the case when $p$ does not divide the order of any maximal subgroup of $M$, and in \cite[Theorem 3A.2; Conjecture 1B.1]{he2025growthproblemsrepresentationsfinite} a formula for $a(n)$ is given in a special case (assuming the `group-injective' condition, see \autoref{def}), with a conjectural answer given in general. The formulas involve character tables and/or decomposition matrices, and are not always easy to apply.   

The contribution of the present paper is as follows. Firstly, we extend the general theory of growth problems, most importantly in \autoref{proj length}, where we generalize \cite[Theorem 3B.2]{he2025growthproblemsrepresentationsfinite} to arbitrary characteristics and hence answer the growth problem for $l(n)$ in full generality. We also discuss how large $n$ needs to be in order for $V^{\otimes n}$ (under certain assumptions) to contain a $kG$-module as submodule or summand, in \autoref{n0 bound} and \autoref{m0 bound}, which implies that eventually almost all factors or summands of $V^{\otimes n}$ are $kG$-modules. We see this 
as an analog of \cite{bryant1972tensor}: when $M$ is a group eventually almost all summands of $V^{\otimes n}$ are projective by \cite{bryant1972tensor}. We also generalize the theory to semigroups in \autoref{semigroup section}.

Secondly, we compute various explicit examples, with a focus on the case where $M$ is one of the diagram monoids, including
\begin{enumerate}
\item The planar diagram monoids: the planar rook monoid $pRo_m$, the Temperley--Lieb monoid $TL_m$, the Motzkin monoid $Mo_m$; and
\item The symmetric diagram monoids: the rook monoid $Ro_m$, the Brauer monoid $Br_m$, the rook Brauer monoid $RoBr_m$, and the partition monoid $Pa_m$. 
\end{enumerate}

\begin{Remark}\label{R:Diagrams}
The diagrammatic representation of these monoids are explained in \cite[Figure (1E.2)]{KhSiTu-monoidal-cryptography}, reproduced here with 
$m$ being the number of strands: 
\begin{gather*}
\begin{tabular}{c|c|c||c|c|c}
\arrayrulecolor{tomato}
Symbol & Diagrams & Name: ${}_{-}$ monoid
& Symbol & Diagrams & Name: ${}_{-}$ monoid
\\
\hline
\hline
$pPa_m$ & \begin{tikzpicture}[anchorbase]
\draw[usual] (0.5,0) to[out=90,in=180] (1.25,0.45) to[out=0,in=90] (2,0);
\draw[usual] (0.5,0) to[out=90,in=180] (1,0.35) to[out=0,in=90] (1.5,0);
\draw[usual] (0.5,1) to[out=270,in=180] (1,0.55) to[out=0,in=270] (1.5,1);
\draw[usual] (1.5,1) to[out=270,in=180] (2,0.55) to[out=0,in=270] (2.5,1);
\draw[usual] (0,0) to (0,1);
\draw[usual] (2.5,0) to (2.5,1);
\draw[usual,dot] (1,0) to (1,0.2);
\draw[usual,dot] (1,1) to (1,0.8);
\draw[usual,dot] (2,1) to (2,0.8);
\end{tikzpicture} & Planar partition
& $Pa_m$ & \begin{tikzpicture}[anchorbase]
\draw[usual] (0.5,0) to[out=90,in=180] (1.25,0.45) to[out=0,in=90] (2,0);
\draw[usual] (0.5,0) to[out=90,in=180] (1,0.35) to[out=0,in=90] (1.5,0);
\draw[usual] (0,1) to[out=270,in=180] (0.75,0.55) to[out=0,in=270] (1.5,1);
\draw[usual] (1.5,1) to[out=270,in=180] (2,0.55) to[out=0,in=270] (2.5,1);
\draw[usual] (0,0) to (0.5,1);
\draw[usual] (1,0) to (1,1);
\draw[usual] (2.5,0) to (2.5,1);
\draw[usual,dot] (2,1) to (2,0.8);
\end{tikzpicture} & Partition
\\
\hline
$Mo_m$ & \begin{tikzpicture}[anchorbase]
\draw[usual] (0.5,0) to[out=90,in=180] (1.25,0.5) to[out=0,in=90] (2,0);
\draw[usual] (1,0) to[out=90,in=180] (1.25,0.25) to[out=0,in=90] (1.5,0);
\draw[usual] (2,1) to[out=270,in=180] (2.25,0.75) to[out=0,in=270] (2.5,1);
\draw[usual] (0,0) to (1,1);
\draw[usual,dot] (2.5,0) to (2.5,0.2);
\draw[usual,dot] (0,1) to (0,0.8);
\draw[usual,dot] (0.5,1) to (0.5,0.8);
\draw[usual,dot] (1.5,1) to (1.5,0.8);
\end{tikzpicture} & Motzkin
& $RoBr_m$ & \begin{tikzpicture}[anchorbase]
\draw[usual] (1,0) to[out=90,in=180] (1.25,0.25) to[out=0,in=90] (1.5,0);
\draw[usual] (1,1) to[out=270,in=180] (1.75,0.55) to[out=0,in=270] (2.5,1);
\draw[usual] (0,0) to (0.5,1);
\draw[usual] (2.5,0) to (2,1);
\draw[usual,dot] (0.5,0) to (0.5,0.2);
\draw[usual,dot] (2,0) to (2,0.2);
\draw[usual,dot] (0,1) to (0,0.8);
\draw[usual,dot] (1.5,1) to (1.5,0.8);
\end{tikzpicture} & Rook Brauer
\\
\hline
$TL_m$ & \begin{tikzpicture}[anchorbase]
\draw[usual] (0.5,0) to[out=90,in=180] (1.25,0.5) to[out=0,in=90] (2,0);
\draw[usual] (1,0) to[out=90,in=180] (1.25,0.25) to[out=0,in=90] (1.5,0);
\draw[usual] (0,1) to[out=270,in=180] (0.25,0.75) to[out=0,in=270] (0.5,1);
\draw[usual] (2,1) to[out=270,in=180] (2.25,0.75) to[out=0,in=270] (2.5,1);
\draw[usual] (0,0) to (1,1);
\draw[usual] (2.5,0) to (1.5,1);
\end{tikzpicture} & Temperley--Lieb
& $Br_m$ & \begin{tikzpicture}[anchorbase]
\draw[usual] (0.5,0) to[out=90,in=180] (1.25,0.45) to[out=0,in=90] (2,0);
\draw[usual] (1,0) to[out=90,in=180] (1.25,0.25) to[out=0,in=90] (1.5,0);
\draw[usual] (0,1) to[out=270,in=180] (0.75,0.55) to[out=0,in=270] (1.5,1);
\draw[usual] (1,1) to[out=270,in=180] (1.75,0.55) to[out=0,in=270] (2.5,1);
\draw[usual] (0,0) to (0.5,1);
\draw[usual] (2.5,0) to (2,1);
\end{tikzpicture} & Brauer
\\
\hline
$pRo_m$ & \begin{tikzpicture}[anchorbase]
\draw[usual] (0,0) to (0.5,1);
\draw[usual] (0.5,0) to (1,1);
\draw[usual] (2,0) to (1.5,1);
\draw[usual] (2.5,0) to (2.5,1);
\draw[usual,dot] (1,0) to (1,0.2);
\draw[usual,dot] (1.5,0) to (1.5,0.2);
\draw[usual,dot] (0,1) to (0,0.8);
\draw[usual,dot] (2,1) to (2,0.8);
\end{tikzpicture} & Planar rook
& $Ro_m$ & \begin{tikzpicture}[anchorbase]
\draw[usual] (0,0) to (1,1);
\draw[usual] (0.5,0) to (0,1);
\draw[usual] (2,0) to (2,1);
\draw[usual] (2.5,0) to (0.5,1);
\draw[usual,dot] (1,0) to (1,0.2);
\draw[usual,dot] (1.5,0) to (1.5,0.2);
\draw[usual,dot] (1.5,1) to (1.5,0.8);
\draw[usual,dot] (2.5,1) to (2.5,0.8);
\end{tikzpicture} & Rook
\\
\hline
$pS_m$ & \begin{tikzpicture}[anchorbase]
\draw[usual] (0,0) to (0,1);
\draw[usual] (0.5,0) to (0.5,1);
\draw[usual] (1,0) to (1,1);
\draw[usual] (1.5,0) to (1.5,1);
\draw[usual] (2,0) to (2,1);
\draw[usual] (2.5,0) to (2.5,1);
\end{tikzpicture} & Planar symmetric
& $S_m$ & \begin{tikzpicture}[anchorbase]
\draw[usual] (0,0) to (1,1);
\draw[usual] (0.5,0) to (0,1);
\draw[usual] (1,0) to (1.5,1);
\draw[usual] (1.5,0) to (2.5,1);
\draw[usual] (2,0) to (2,1);
\draw[usual] (2.5,0) to (0.5,1);
\end{tikzpicture} & Symmetric
\\
\end{tabular}
.
\end{gather*}
The left half of the table above contains \emph{planar} monoids, the right half 
\emph{symmetric} monoids. These monoids appear under various different names in the literature, see \cite[Remark 1E.1]{KhSiTu-monoidal-cryptography} for some details.

The planar monoids have trivial group of units, while the symmetric monoids have the symmetric groups as the group of units. We note that the planar partition monoid $pPa_m$ is isomorphic to $TL_{2m}$, and the two bottom monoids are groups, and all these are therefore omitted.
\end{Remark}

\begin{Remark}\label{R:Diagrams2}
We have two more symmetric diagram monoids that we study:
\begin{gather*}
T_m\colon
\begin{tikzpicture}[anchorbase,scale=0.55]
\draw[usual] (4,-3) to[out=90,in=270] (1,-1);
\draw[usual] (3,-3) to[out=90,in=270] (1,-1);
\draw[usual] (2,-3) to[out=90,in=270] (1,-1);
\draw[usual] (1,-3) to[out=90,in=270] (1,-1);
\draw[usual] (0,-3) to[out=90,in=270] (-1,-1);
\draw[usual] (-1,-3) to[out=90,in=270] (-3,-1);
\draw[usual] (-2,-3) to[out=90,in=270] (-1,-1);
\draw[usual] (-3,-3) to[out=90,in=270] (-2,-1);
\draw[usual,dot] (0,-1) to (0,-1.25);
\draw[usual,dot] (2,-1) to (2,-1.25);
\draw[usual,dot] (3,-1) to (3,-1.25);
\draw[usual,dot] (4,-1) to (4,-1.25);
\end{tikzpicture}
,\quad
PT_m\colon
\begin{tikzpicture}[anchorbase,scale=0.55]
\draw[usual] (4,-3) to[out=90,in=270] (1,-1);
\draw[usual] (3,-3) to[out=90,in=270] (1,-1);
\draw[usual] (2,-3) to[out=90,in=270] (1,-1);
\draw[usual] (1,-3) to[out=90,in=270] (1,-1);
\draw[usual] (0,-3) to[out=90,in=270] (-1,-1);
\draw[usual] (-1,-3) to[out=90,in=270] (-3,-1);
\draw[usual] (-2,-3) to[out=90,in=270] (-1,-1);
\draw[usual,dot] (-2,-1) to (-2,-1.25);
\draw[usual,dot] (-3,-3) to (-3,-2.75);
\draw[usual,dot] (0,-1) to (0,-1.25);
\draw[usual,dot] (2,-1) to (2,-1.25);
\draw[usual,dot] (3,-1) to (3,-1.25);
\draw[usual,dot] (4,-1) to (4,-1.25);
\end{tikzpicture}
,
\end{gather*}
generated by crossings, merges, top dots, for $T_m$, and additionally bottom dots for $PT_m$; see e.g. 
\cite[Section 4]{Tu-sandwich-cellular} for details. These are called 
full transformation monoid $T_m$ and partial transformation monoid $PT_m$, respectively.
\end{Remark}

In \autoref{diagram bn}, we give a formula for $a(n)$ when $k=\C$ and $M$ is any of the diagram monoids listed above. For the planar monoids, we in addition obtain explicit formulas for $l(n)$ (in fact for $[V^{\otimes n}: S]$ where $S$ is any simple representation) in \autoref{planar}. These are nontrivial, and along the way we obtain the monoid character tables of the planar monoids and their inverses.
To the best of our knowledge, this paper is the first systematic study of tensor products of representations of the above diagram monoids. 

We also establish that many important classes of monoids satisfy \cite[Conjecture 1B.1]{he2025growthproblemsrepresentationsfinite}, by showing that they satisfy the group-injective condition or otherwise. These include the diagram monoids in \autoref{R:Diagrams} and \autoref{R:Diagrams2} as well as inverse monoids, the full linear monoid $\M_m(q)$, and the Catalan monoid $Ca_m$, and some other monoids.

In summary, we have a formula for $l(n)$ in all cases, and for $b(n)$ for the following cases. Here groups or monoids are always finite, representations are always finite dimensional, and the references indicate how to obtain a formula.
\begin{enumerate}

\item All finite groups. 
\cite{he2024growthproblemsrepresentationsfinite}.

\item All finite semisimple monoids. \cite{he2025growthproblemsrepresentationsfinite}.

\item The above diagram monoids in characteristic 0. Mostly this paper.

\item All inverse monoids. This paper.

\item Specific example such as all monoids up to order $\leq 7$ and several other small monoids. Also Catalan monoids, full linear monoids, (submonoids of) $PT_m$. \cite{he2025growthproblemsrepresentationsfinite} and this paper.

\end{enumerate}
There is also work for infinite groups, which is mostly in characteristic zero, see e.g. \cite{Bi-asymptotic-lie}, or in specific cases in finite characteristic, see e.g. \cite{coulembier2024fractalbehaviortensorpowers,larsen2024boundsmathrmsl2indecomposablestensorpowers}.

\begin{Remark}
We have used Magma (cf. \cite{bosma1997magma}) for the computation of particular examples. See \cite{code} for details of code implementations.  
\end{Remark}

\begin{Remark}
This paper is fairly self-contained, but we assume some basics from monoid representation theory (as for example in \cite{steinberg2016representation}) and some background on diagram monoids (as e.g. in \cite{KhSiTu-monoidal-cryptography}).
\end{Remark}

\noindent\textbf{Acknowledgments.} We would like to express our gratitude to Volodymyr (Walter) Mazorchuk for insightful discussions, in particular for suggesting that we look at the partial transformation monoid. DT was supported by the ARC Future Fellowship FT230100489 and the persistent illusion that diagram monoids are fun.

\section{Growth rate of length of tensor powers}
\subsection{General theory}
Let $\{e_i\mid 1\le i \le m\}$ denote a complete set of idempotent representatives for the regular $\mathcal{J}$-classes of $M$, with corresponding maximal subgroups $G_{e_i}$. In the case when $\text{char\ }k =p \nmid \lvert G_{e_i}\rvert$ for each $1\le i\le m$ (so each $kG_{e_i}$ is semisimple), an exact formula for $l(n)$ is known from \cite[Theorem 3B.2]{he2025growthproblemsrepresentationsfinite}. We now extend this formula to arbitrary characteristic.

Let $L$ denote the \textit{decomposition matrix} of $M$ corresponding to the restriction isomorphism of Grothendieck rings $\Res: G_0(kM)\to \prod_{i=1}^m G_0(kG_{e_i}), [V] \mapsto ([e_1V],\dots,[e_mV])$, and let $Y$ denote the block triangular matrix \[
Y = \begin{pmatrix}
Y_1 & 0      & \cdots & 0      \\[6pt]
0      & Y_2 & \cdots & 0      \\[6pt]
\vdots & \vdots & \ddots & \vdots \\[6pt]
0      & 0      & \cdots & Y_m
\end{pmatrix}
\]
where for $1\le i\le m$, $Y_i$ is the projective (Brauer) character table of $kG_{e_i}$. We observe that when $\text{char\ }k \nmid \lvert G_{e_i}\rvert$ for each $i$, the projective character tables are just the ordinary character tables, and $Y$ is the matrix $X$ such that $L^{T}X$ is the (ordinary) character table of $M$ (see \cite[Corollary 7.17]{steinberg2016representation}). Assume that for $1\le i \le m$, $G_{e_i}$ has $k_i$ ($p$-regular) conjugacy classes, which we label by $C_{i,j}$, $1\le j \le k_i$, with representatives $g_{i,j}$. Let $V_{i,j}$ denote the simple $kM$-module indexed by $C_{i,j}$, and let $Y_{i,j}$ denote the column of $Y$ corresponding to $C_{i,j}$, regarded as a column vector. If $V$ is a $kM$-module, let $Z_V(M)$ denote the set of $x\in M$ acting as a scalar on $V$, and let $\omega_V(x), x\in Z_V(M)$, denote the corresponding scalar. If $V$ has character $\chi$, let $\chi_{\textup{sec}}$ be $0$ if $Z_V(M)=M$ and otherwise let it denote any second largest character value (in terms of modulus) that $\chi$ takes on (the largest being the dimension). Finally, we refer to \cite[\S 5]{steinberg2023modularrepresentationtheorymonoids} for the definition of Brauer characters for monoid representations.

\begin{Theorem}\label{proj length}
Let $V$ be a $kM$-module with (Brauer) character $\chi$, then:
\begin{enumerate}
\item For each simple $kM$-module $V_{r,s}$, we have \begin{equation}\label{eqn: individual term}
[V^{\otimes n}: V_{r,s}] =\sum_{i=1}^m \frac{1}{\lvert G_{e_i}\rvert} \sum_{j=1}^{k_i}\lvert C_{i,j} \rvert T^{r,s}_{i,j} (\chi(g_{i,j}))^n,
\end{equation} where $T^{r,s}_{i,j}$ denotes the $(r,s)$th entry of $L^{-1}\overline{Y}_{i,j}$ (here the bar denotes complex conjugation), and so
\begin{equation}\label{eqn:exact}
l(n)=\sum_{i=1}^m \frac{1}{\lvert G_{e_i}\rvert }\sum_{j=1}^{k_i}\lvert C_{i,j}\rvert T_{i,j} \big( \chi(g_{i,j}) \big)^n         
\end{equation}
where $T_{i,j}$ is sum over the entries of ${L^{-1}}\overline{Y_{i,j}}$.     
\item We have the asymptotic formula \begin{equation}\label{eqn:asym}
l(n)\sim k(n):=\sum_{i=1}^m \frac{(\dim V)^n}{\lvert G_{e_i}\rvert }\sum_{j: g_{i,j} \in Z_V(G_{e_i})}\lvert C_{i,j}\rvert T_{i,j} \big( \omega_V(g_{i,j}) \big)^n.       
\end{equation} 
An asymptotic formula for $[V^{\otimes n}: V_{r,s}]$ is obtained by replacing $T_{i,j}$ everywhere with $T^{r,s}_{i,j}$.
\item $\lvert l(n)/k(n)-1\rvert \in \mathcal{O}(\lvert \chi_{\mathrm{sec}}/\dim V \rvert^n)$, and
\item $\lvert l(n)-k(n) \rvert \in \mathcal{O}(\lvert \chi_{\mathrm{sec}}\rvert^n)$.
\end{enumerate}
\end{Theorem}

\begin{proof}
First assume $M=G$. Let $J=J(G)$ be the matrix $J_{ij}=[V\otimes V_j: V_i]$ where $1\le i,j\le N$ index the simple $kG$-modules. The matrix $J$ represents the linear map $\chi_W\mapsto \chi_{V\otimes W}$ in the (Brauer) character ring, with respect the basis of simple characters; by writing the same map in the basis of the projective (Brauer) characters, we obtain $J=CKC^{-1}$ where $C$ is the Cartan matrix and $K$ is the matrix such that $K_{ij}$ is the number of copies of $P_i$ showing up as a direct summand of $V\otimes P_j$. (The $P_i$'s for $1\le i \le N$ are the projective indecomposable $kG$-modules.) It is shown in \cite[Lemma 7]{he2024growthproblemsrepresentationsfinite}, and easy to verify directly, that for $1\le i\le N$, $\chi(g_i)$ is an eigenvalue of $K$ with left eigenvector $\begin{pmatrix}\chi_{P_1}(g_i)& \dots \chi_{P_N}(g_i)\end{pmatrix}$ and right eigenvector  ${\begin{pmatrix}\overline{\chi_{V_1}(g_i)} & \dots \overline{\chi_{V_N}(g_i)}\end{pmatrix}}$, where $g_i$ is a representative of the $i$th conjugacy class. The left (resp. right) eigenvectors of $J$ are obtained from the left (resp. right) eigenvectors of $K$ by multiplying $C^{-1}$ (resp. $C$) on the right (resp. left), and hence are $w_i^T:=\begin{pmatrix}
\chi_{V_1}(g_i) & \dots & \chi_{V_N}(g_i)    
\end{pmatrix}$ and  $v_i:=\begin{pmatrix}
\overline{\chi_{P_1}(g_i)} & \dots & \overline{\chi_{P_N}(g_i)}    
\end{pmatrix}$, respectively, with corresponding eigenvalue $\chi_V(g_i)$ and for $1\le i \le N$. Observe that the $v_i$'s are complex conjugates of columns of the projective (Brauer) character table of $G$.

Now consider the matrix $J(M)$ for the monoid defined analogously but now with rows and columns indexed by simple $kM$-modules $V_{i,j}$. By using $L$ to go between $J(M)$ and the $J(G_{e_i})$'s as in \cite[Lemma 3B.1]{he2025growthproblemsrepresentationsfinite}, we check that the right eigenvectors of $J(M)$ are exactly the vector $\overline{Y_{i,j}}$'s, whereas the left eigenvectors are of the form $$w^T_{i,j}:=\begin{pmatrix} \chi_{V_{1,1}}(g_{i,j}) & \dots & \chi_{V_{m,k_m}}(g_{i,j}) \end{pmatrix},$$ with corresponding eigenvalue $\chi(g_{i,j})$. We thus obtain the eigendecomposition
\begin{equation}\label{eigendecomp}(J(M))^n=\sum_{i=1}^m \frac{1}{\lvert G_{e_i}\rvert}\sum_{j=1}^{k_i}\lvert C_{i,j}\rvert w^T_{i,j}\overline{Y_{i,j}} (\chi(g_{i,j}))^n,
\end{equation}
where the factors $\lvert C_{i,j}\rvert/\lvert G_{e_i}\rvert$ ensure that the matrices in \autoref{eigendecomp} are projections. 
Assume without loss of generality that $V_{1,1}$ is the trivial $kM$-module, then by usual Perron--Frobenius theory (see \cite[Theorem 5.10]{lacabanne2024asymptotics}) $[V^{\otimes n}: V_{i,j}]$ is exactly the sum over all the matrix summands in \autoref{eigendecomp} of the entry whose column index is $(1,1)$ and whose row index is $(i,j)$. This implies item a) of the Theorem. Item b) follows by observing that the eigenvalues (character values) with maximal modulus ($=\dim V$) will have their $n$th power dominating the expression for $b(n)$ when $n$ is sufficiently large, and items c) and d) follow by considering the expressions $b(n)-a(n)$ and $b(n)/a(n)$. 
\end{proof}

If no element apart from $1$ acts as identity, formula \autoref{eqn:asym} yields \begin{equation} \label{eqn:no identity length formula}
l(n) \sim k(n)=\frac{1}{\lvert G\rvert}\sum_{\substack{1\le t\le N \\ g_t \in Z_V(G)}} T_{t}\big(\omega_V(g_t)\big)^n\cdot (\dim V)^n,    
\end{equation}
where $1\le t\le N$ index the conjugacy classes of the group of units $G$, and $T_t$ is the sum over the conjugate of the $t$th column of the projective (Brauer) character table of $G$. (Note that the submatrix of $L$ corresponding to $G$ is the identity matrix, so in this case $L$ no longer figures in the formula. Also the $\lvert C_t\rvert$'s disappear since the classes are central). In other words, in this case the asymptotic growth rate of $l(n)$ is the same as that of $l(n)$ for $V$ considered as a $kG$-module by restriction. We note also that if $M=G$, \autoref{eqn:exact} simplifies to \begin{equation}\label{eqn:group exact length}
l(n) =\frac{1}{\lvert G\rvert}\sum_{1\le t\le N } \lvert C_t\rvert T_{t}\big(\chi_V(g_t)\big)^n.   
\end{equation}

\begin{Remark}
The exact formula in \autoref{proj length} generalizes \cite[Theorem 3B.2]{he2025growthproblemsrepresentationsfinite} (for monoids in good characteristics), which in turn generalizes \cite[Theorem 3]{he2024growthproblemsrepresentationsfinite} (for groups in nonmodular characteristics). We also extend previous results concerning the asymptotic formula $k(n)\sim l(n)$ in the case when $M=G$ obtained in \cite[Proposition 2A.2]{coulembier2024fractalbehaviortensorpowers}, \cite[(2A.1)]{lacabanne2023asymptotics} and \cite[Proposition 2.1]{coulembier2023asymptotic}.
\end{Remark}

\begin{Proposition}
Let $V$ be a $kM$-module on which no element apart from $1\in M$ acts as identity, then:
\begin{enumerate}
\item If $Z(G)$ is either trivial or a $p$-group (where $p>0$ is the field characteristic), then  \begin{equation}\label{eqn: no period length}
l(n)\sim k(n)=\frac{1}{\lvert G\rvert} \sum_{i=1}^N \dim P_i \cdot (\dim V)^n,   
\end{equation}
where $P_i, 1\le i \le N$ are the projective indecomposable $kG$-modules.
\item We have \begin{equation}\label{eqn: trivial asym length}
l(n)\sim k(n)=(\dim V)^n
\end{equation}
if and only if either $p\nmid \lvert G\rvert$ and $G$ is abelian or $p\mid \lvert G\rvert$ and $G$ is the extension of a p-group by an abelian group.  
\end{enumerate}
\end{Proposition}

\begin{proof}
The first statement follows from \autoref{eqn:no identity length formula} because $Z_V(M)=\{1\}$ by hypothesis. To see that the second statement is true, note that $\sum_{i=1}^N \dim P_i = \lvert G\rvert$ if and only if all simple $kG$-modules are one-dimensional. This happens in nonmodular characteristic if and only if $G$ is abelian, and happens in positive characteristic if and only if $G$ is the extension of a $p$-group by an abelian group. This last fact is well-known, see e.g.\cite[Theorem 4.3]{almeida2009representation} which gives a more general result characterizing semigroups all of whose simple representations are one-dimensional. 
\end{proof}
\begin{Remark}
If all simple $kM$-modules are one-dimensional (that is, $kM$ is a split basic algebra), then we even have $l(n)=(\dim V)^n$ exactly. The monoids (in fact semigroups) satisfying this property are characterized in \cite[Theorem 4.3]{almeida2009representation}. 
\end{Remark}
\begin{Remark}
The formula \autoref{eqn: trivial asym length} applies to the planar diagram monoids, whose group of units are trivial. In \autoref{planar} we will study the harder problem of determining an exact formula for $l(n)$ as well as quantifying the rate of convergence and variance.
\end{Remark}

Recall that the fusion graph $\Gamma_l$ (for length) associated with $V$ is the directed and weighted graph whose vertices are the simple $kM$-modules, such that we draw an edge of weight $m$ from $W$ to $X$ if $[V\otimes W: X]=m$. By \autoref{eqn:no identity length formula}, if we assume no element apart from $1\in M$ acts as identity on $V$, then the asymptotic growth rate $k(n)\sim l(n)$ is completely determined by the group $G$. This corresponds to the fact that the simple $kG$-modules (which are simple $kM$-modules by letting non-units act as 0) form a strongly-connected component $\Gamma_l^G$ of $\Gamma_l$ which is reachable from any vertex, and which has no edges leaving the component; for $n$ large, this component determines $l(n)$ (see \cite[Proposition 4.22]{lacabanne2024asymptotics}). We note that $\Gamma_l^G$ is strongly-connected because it is isomorphic to the fusion graph of $V$ as a (faithful!) $kG$-module. We can bound the number of steps needed to go from the trivial $kM$-module to $\Gamma_l^G$, which may be seen as an alternative measure of the rate of convergence $k(n)\sim l(n)$.

\begin{Proposition}\label{n0 bound}
Let $V$ be a $kM$-module on which no element apart from $1\in M$ acts as identity. Let $n_0$ be the smallest natural number such that $[V^{\otimes n_0}:V_i]>0$ for some $V_i$ which is (the induction of) a simple $kG$-module, then $n_0 \le L-1$ where $L$ is the number of $\mathcal{L}$-classes of $M$. 
\end{Proposition}
\begin{proof}
Let $Z_n$ denote the submodule of $V^{\otimes n}$ consisting of all vectors on which $M\setminus G$ acts as 0. If $m_1,\dots, m_{L-1}$ are representatives for the $L-1$ $\mathcal{L}$-classes which are not $G$, and $v_i\neq 0$ is in the null space of the action of $m_i$, for $1\le i \le L-1$, then the vector $v_1\otimes \dots \otimes v_{L-1}\neq 0$ is in $Z_{L-1}$, since we can write any $m\in M\setminus G$ as $m=m'm_i$ for some $m'\in M\setminus G$, $1\le i \le L-1$, and so $m$ acts as 0 on this vector. It follows that $V^{\otimes (L-1)}$ has a composition factor (in fact a submodule) which is a $kG$-module.  
\end{proof}

We now consider some explicit examples. 

\subsection{Length for planar monoids}\label{planar}
From \autoref{eqn: trivial asym length}, we know that for the planar monoids $pRo_m, TL_m$ and $Mo_m$ the asymptotic formula $k(n)\sim l(n)$ is simply $(\dim V)^n$ (assuming no element apart from 1 acts as identity). Therefore, our focus will be on computing the exact formula $l(n)$, and quantifying the rate of convergence and variance. 

Let $\mathcal{X}_m=pRo_m$, $TL_m$ or $Mo_m$. The $\mathcal{J}$-classes of $\mathcal{X}_m$ consist of all diagrams with the same number of through strands $k$ (see e.g.\cite[\S 4B]{KhSiTu-monoidal-cryptography}). Thus, the maximal subgroups (which are trivial) and $\mathcal{J}$-classes (which are also the generalized conjugacy classes) of $\mathcal{X}_m$ are indexed by an integer $l\in \Lambda_m$. For $\mathcal{X}_m=pRo_m$ or $Mo_m$ we have $\Lambda_m =\{0,1,\dots,m\}$, whereas for $\mathcal{X}_m=TL_m$, $\Lambda$ is either $\{0,2,\dots, m-2, m\}$ or $\{1,3\dots, m-2, m\}$, depending on the parity of $m$. For $l\in \Lambda_m$, let $V_l$ denote the corresponding simple $k\mathcal{X}_m$-modules, with character $\chi_l$, and let $S_l$ denote the corresponding indecomposable Schützenberger representations. We note that $S_l$'s are also the \textit{cell representation} of the monoid algebras $k\mathcal{X}_m$ coming from their cellular algebra structure.  We will denote by $\chi_i(j)$ the value of the character $\chi_i$ at the $\mathcal{J}$-class indexed by $j$. Finally, we observe that the character table of $\mathcal{X}_m$ is upper triangular when both rows and columns are indexed by $l \in \Lambda_l$ with $l$ in the natural order. The matrix $L$ is the transpose of the character table, and $Y$ is the identity matrix.  
\subsubsection{Planar rook monoid}
The planar rook monoid $pRo_m$ has order $\binom{2m}{m}$ and is semisimple over any field $k$. 
\begin{Lemma}
The simple $kPRo_m$ characters are given by $$\chi_i(j)=\begin{cases}
\binom{j}{i} & i\le j,\\
0 & i>j. 
\end{cases}$$
That is, the character table is the Pascal's triangle.
\end{Lemma}
\begin{proof}
The result was stated in \cite[Theorem 5.1]{flath2009planar} for $k=\C$ but (by passing to Brauer characters if necessary) it is actually valid over any field $k$, since the simple representations are all permutation representations. 
\end{proof}
We can now compute $l(n)$.

\begin{Proposition}\label{plannarrook formula}
Let $V$ be a $kpRo_m$-module with character $\chi$. We have:
\begin{enumerate}
\item $$[V^{\otimes n}: V_l]=\sum_{j=0}^l (-1)^{l-j}\binom{l}{j}(\chi(j))^n,$$
and so
$$l(n)=b(n)=\sum_{l=0}^m \sum_{j=0}^l (-1)^{l-j}\binom{l}{j}(\chi(j))^n.$$
\item If $V=V_i$, then $l(n)\sim k(n)=(\dim V)^n=(\binom{m}{i})^n$, $\lvert l(n)/k(n)-1\rvert \in \mathcal{O}((\frac{m-i}{m})^n)$, and $\lvert l(n)-k(n)\rvert \in \mathcal{O}((\binom{m-1}{i})^n)$.
\end{enumerate}
\end{Proposition}

\begin{proof}
From the identity $\sum_{k=0}^m (-1)^{j-k}\binom{j}{k}\binom{k}{i}=\delta_{i,j}$, it can be deduced that $L^{-1}_{ij}=(-1)^{i-j}\binom{i}{j}$. Now the result follows by applying \autoref{proj length} and observing that for $V=V_i$, $\chi_{\textrm{sec}}=\binom{m-1}{i}$. 
\end{proof}

\begin{figure}[H]
\centering
\begin{minipage}{0.5\textwidth}
\centering
\includegraphics[width=0.9\linewidth]{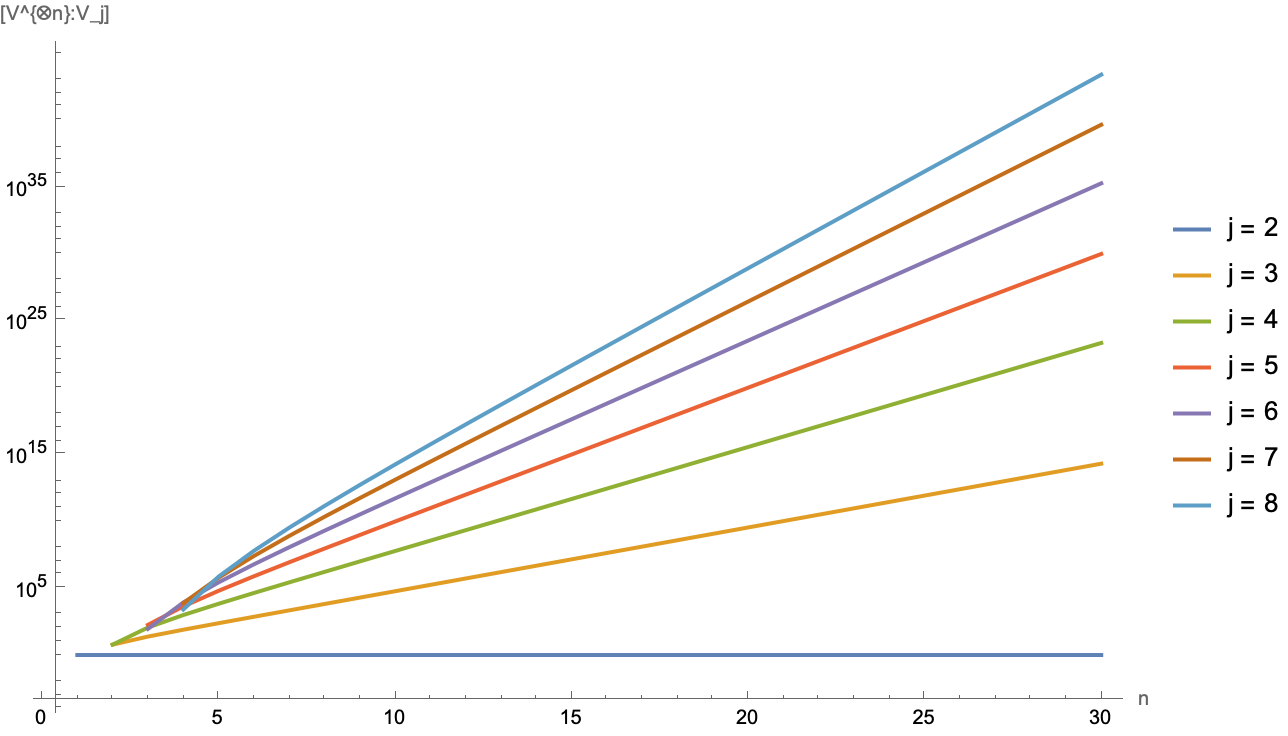}
\end{minipage}%
\begin{minipage}{0.5\textwidth}
\centering
\includegraphics[width=0.9\linewidth]{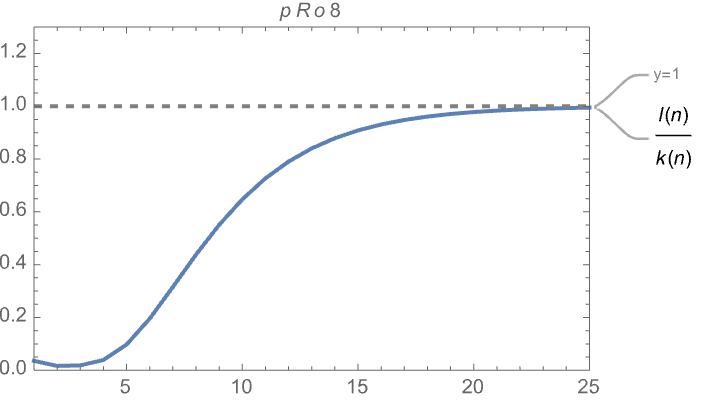}
\end{minipage}
\caption{The growth rate of the composition multiplicities and the ratio $l(n)/k(n)$ for $V=V_2$ for $pRo_8$.}
\label{fig:planarrook}
\end{figure}

\begin{figure}[H]
\centering
\begin{minipage}{0.40\textwidth}
\centering
\includegraphics[width=0.6\linewidth]{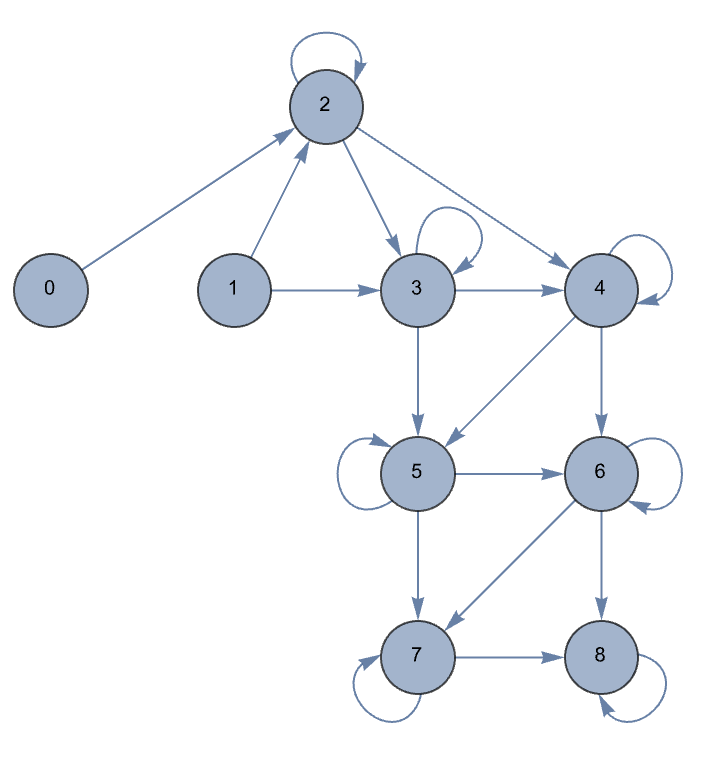}
\end{minipage}
\begin{minipage}{0.4\textwidth}
\centering

\[
\begin{psmallmatrix}
0 & 0 & 0  & 0  & 0  & 0  & 0  & 0  & 0 \\
0 & 0 & 0  & 0  & 0  & 0  & 0  & 0  & 0 \\
1 & 2 & 1  & 0  & 0  & 0  & 0  & 0  & 0 \\
0 & 3 & 6  & 3  & 0  & 0  & 0  & 0  & 0 \\
0 & 0 & 6  & 12 & 6  & 0  & 0  & 0  & 0 \\
0 & 0 & 0  & 10 & 20 & 10 & 0  & 0  & 0 \\
0 & 0 & 0  & 0  & 15 & 30 & 15 & 0  & 0 \\
0 & 0 & 0  & 0  & 0  & 21 & 42 & 21 & 0 \\
0 & 0 & 0  & 0  & 0  & 0  & 28 & 56 & 28
\end{psmallmatrix}
\]

\end{minipage}
\caption{The fusion graph and its adjacency matrix for $V=V_2$ as a $kpRo_8$-module. The shortest path to go from the trivial module $V_0$ to $V_8$ has length $4=\ceil{8/2}$=4.}
\label{fig:planarrookfusion}
\end{figure}

Thus, the convergence for $V=V_i, i>0$ slows as $i$ gets smaller. We can also determine the tensor product decompositions of simple modules, and the value of $n_0$.
\begin{Proposition}\ 
\begin{enumerate}
\item We have $$[V_i\otimes V_j: V_l]=
\binom{l}{i}\binom{i}{i+j-l},$$ in particular $[V_i\otimes V_j:V_l]=0$ if and only if $l<\max(i,j)$ or $l>\min(i+j)$.
\item In the language of \autoref{n0 bound}, for $V=V_i$ we have $n_0=\ceil{m/i}$ where $\ceil{\cdot}$ is the ceiling function. 
\end{enumerate}
\end{Proposition}
\begin{proof}
Item a) follows from the identity, easily verified, that $$\binom{r}{i}\binom{r}{j}=\sum_{l=0}^{i+j}\binom{l}{i}\binom{i}{i+j-l}\binom{r}{l},$$ which is true for each $0\le r \le m$. From this, we see that the largest $l$ such that $V_l$ is a composition factor of $V^{\otimes n}$ is $\min\{ni,m\}$, so item b) follows. 
\end{proof}
We plot in the left panel of \autoref{fig:planarrook} the growth rate of $[V^{\otimes n}:V_j]$ for varying $j$ when $V=V_2$ for $pRo_8$; the vertical axis is on a logarithmic scale. On the right panel we plot the corresponding ratio $l(n)/k(n)$, which converges to 1, as expected.  We also plot the fusion graph (with labels suppressed) and its adjacency matrix in \autoref{fig:planarrookfusion}; note that $n_0=\ceil{m/i}=4$ in this case. 

\subsubsection{Temperley--Lieb monoid}\label{SS:TL}
Let $D_{ij}=[S_j:V_i]$ denote the (cellular) decomposition matrix encoding the multiplicities of the simple $kTL_n$-modules in the cell modules. Over a field of arbitrary characteristic, $D$ and its inverse are known from \cite[Theorem 8.4; Theorem 9.2]{Spencer_2023}. (For clarity, we are always using the $l=3$ case of Spencer's results.) Thanks to this, to invert the character table $X$ of $TL_n$ it suffices to invert the `cell character table', whose rows are the characters of the $S_i$'s. We note that the cell modules and their characters do not depend on the field characteristic. 
\begin{Proposition}\label{standard characters TL}
Let $\chi_{S_i}$ denote the character of $S_i$, and write $c(j,i)$ for $\frac{j-i}{2}$. We have 
$$\chi_{S_i}(j)=\alpha_{j,i}:=\begin{cases} \frac{j-2c(j,i)+1}{j-c(j,i)+1}\binom{j}{c(j,i)} & j \ge i, \\
0 & else.
\end{cases}$$
\end{Proposition}
\begin{proof}
These characters are the same as the simple characters in the semisimple setting, which were computed in \cite[\S 2]{halverson1995characters}. To quickly verify, note that the cell modules have a basis given by half diagrams with $m$ points and $i$ through strands (`$(m,i)$ half diagrams'). The $(m,i)$ half diagrams fixed by an element in $TL_m$ with $j$ through stands are just the $(j,i)$ half diagrams with the appropriate number of cups added. The number of $(j,i)$ half diagrams is the dimension of $S_i$ as a $kTL_j$-module, which is just $\alpha_{j,i}$ (see e.g.\cite[Proposition 4B.3]{KhSiTu-monoidal-cryptography}). 
\end{proof}

\begin{Remark}\label{tl remark}
The formula in \autoref{standard characters TL} is independent of $m$, and so the cell character table of $TL_m$ naturally embed into that of $TL_m'$, where $m'>m$ and $m$ and $m'$ have the same parity. Thus, the cell character table of any $TL_m$ is the truncation of some infinite uni upper triangular matrix $S$, with the $i$th row given by $\alpha_{j,i}$. In fact, if $i$ is even (resp. odd), the row of $S$ indexed by by $i$ corresponds to the $(i/2)$th (resp. $((i-1)/2)$th) convolution of the Catalan numbers, which are the falling diagonals of the Catalan triangle, cf. \cite[A009766]{Oeis}. 
\end{Remark}

Let $S$ denote the infinite matrix defined in \autoref{tl remark} (for either parity). If $F$ is a generating function in $x$, we write $[x^n]F$ for the coefficient of $F$ in front of $x^n$.
\begin{Lemma}\label{inverse TL}
The inverse of $S$ is $$S^{-1}_{ij}=\begin{cases}
[x^{(j-i)/2}](1+x)^{-(i+1)}=(-1)^{\frac{j-i}{2}}\binom{\frac{i+j}{2}}{i} & j\ge i, \\
0 & j< i.
\end{cases}$$    
\end{Lemma}
\begin{proof}
The observation in \autoref{tl remark} that the rows of $S$ are convolutions of the Catalan numbers can be rephrased by saying that $S$ is the Riordan array $(C(x),xC(x))$ where $C(x)=1+x+2x^2+5x^3+\dots$ is the generating function for the Catalan numbers, satisfying $C(x)=1+x(C(x))^2$. It follows that $S^{-1}$ is exactly the inverse Riordan array $(1/(1+x),x/(1+x))$, whence the claim follows by a computation. (We refer to \cite{SHAPIRO1991229} for details on the theory of Riordan groups.)
\end{proof}

Let $k$ be a field of characteristic $p$, and let $q\in k^\times$ be such that $q+q^{-1}=1$. (For the purpose of this section, set $p=\infty$ if $k$ has characteristic zero.) For $a\in \N$, the quantum integer $[a]_q$ is defined as $[0]_q=0,[a]_q=q^{-(a-1)}+q^{-(a-3)}+\dots+q^{(a-3)}+q^{(a-1)}$. Write $l \in \Z_{\ge 0}\cup \{\infty\}$ for the \textit{quantum characteristic} of $TL_m$ over $k$ (that is, the least $n$ such that the quantum integer $[n]_q$ vanishes; set $l=\infty$ if the quantum integers never vanish.) For $a\in \N$, we have the $(l,p)$-expansion $$a=[a_t,a_{t-1},\dots,a_0]_{p,l}=\sum_{i=0}^t a_ip^{(i)},$$
where $p^{(i)}=lp^{i-1}$ with the understanding $p^{(0)}=1$. We require $0\le a_0 < l$, and $0\le a_i<p$ for $i>0$. If $a+1=[a_t,a_{t-1},\dots,a_0]_{p,l}$, the support of $a$ is defined as $\textrm{supp}(a)=\{a_tp^{(t)}\pm a_{t-1}p^{(t-1)}\pm \dots \pm a_0p^{(0)}-1\}$. (We only need this for $l=3$ but its somehow nicer to write this out for general $l$.)
\begin{Proposition}\label{TL length}
Let $V$ be a $kTL_m$-module with (Brauer) character $\chi$, then $$[V^{\otimes n}: V_i]=\sum_{z \in \textup{supp}(i)} \sum_{\substack{0 \le j \le z \\ j \equiv z \pmod{2}}}
\bigl(\chi(j)\bigr)^{n}\,
(-1)^{\frac{z-j}{2}}\,
\binom{\frac{j+z}{2}}{j},$$ and so
\[
l(n)
=
\sum\limits_{\substack{0\le i\le m\\i\equiv m\pmod2}}
\;\sum\limits_{z\in\mathrm{supp}(i)}
\;\sum\limits_{\substack{0\le j\le z\\j\equiv z\pmod2}}
(\chi(j))^n\,
(-1)^{\tfrac{z-j}{2}}\,
\binom{\tfrac{j+z}{2}}{j}.
\]
\end{Proposition}
\begin{proof}
Using \autoref{inverse TL} we may write $\chi$ in the basis of characters of cell modules, so that $$\chi = \sum_{i} c_i \chi_{S_i}, c_i=\sum_{j} \chi(j)(-1)^{\frac{j-i}{2}}\binom{\frac{i+j}{2}}{i}.$$ By \cite[Theorem 8.4] {Spencer_2023}, $V_{i}$ is a composition factor of $S_{z}$ if and only if $z\in \mathrm{supp}(i)$, in which case it appears with multiplicity 1. Using this, we can write the characters of cell modules in the basis of the characters of simple modules, and the claims follow. 
\end{proof}

In characteristic zero, we have $(p,l)=(\infty,3)$ and the $(p,l)$ expansion has at most two terms. It follows that $\mathrm{supp}(i)$ is either $\{i\}$ or $\{i,i^-\}$ (on the number line, $i^-$ corresponds to the leftward reflection of $i$ across the nearest `critical integer', i.e. an integer which is 2 mod 3; we define $i^+$ analogously, by rightward reflection). Let $P_i$ denote the projective cover of $V_i$. As explained in e.g.\cite[\S 2]{Bellet_te_2017}, we actually have the short exact sequences \begin{equation}\label{eqn: SES ir TL}
0 \to V_{i^+} \xrightarrow{} S_i \xrightarrow{} V_i \to 0
,\end{equation} and
\begin{equation}\label{eqn:SES proj TL}
0 \to S_{i^-} \xrightarrow{} P_i \xrightarrow{} S_i \to 0
,\end{equation}
where $V_i^{+}$ (resp. $V_i^-$) is potentially zero, which happens if either $j$ is critical or if $i^{+}$ (resp. $i^{-}$) is out of range, i.e is not in $\Lambda_m$. This allows us to easily determine the characters of the simple and projective modules. Write $i^{(0)}$ for the largest  non-critical integer in $\Lambda_m$, and denote the successive leftward reflections of $i^{(0)}$ by $i^{(1)}, i^{(2)}, $ etc. As before, write $\chi_{i}$ for the character of $V_i$, and write $\phi_i$ for the character of $P_i$. 

\begin{Proposition}\label{TL irr and proj char}
Let $k$ have characteristic zero. If $i$ is critical, then $\phi_i=\chi_i=\chi_{S_i}$. If $i=i^{(t)}$ is non-critical, we have \begin{enumerate}
\item $\chi_{i^{(t)}}(r)=\sum_{s=0}^t (-1)^s \alpha_{r,i^{(t-s)}},$ and
\item $\phi_i^{(t)}(r)=\alpha_{r,i^{(t)}}+\alpha_{r,i^{(t+1)}}$,  
\end{enumerate} 
where $\alpha_{r,i^{(j)}}=0$ if $i^{(j)}\notin \Lambda_m$.  
\end{Proposition}\label{TL char}
\begin{proof}
Item (b) follows from equation \autoref{eqn:SES proj TL}, and item (a) follows from repeatedly applying \autoref{eqn: SES ir TL}: we know $\chi_{S_{i^{(0)}}}=\chi_i^{(0)}$ and $\chi_{i^{(t+1)}}=\chi_{S_{i^{(t)}}}-\chi_{i^{(t)}}$.  
\end{proof}

\begin{Remark}
It is not difficult to generalize \autoref{TL irr and proj char} to arbitrary characteristic, since we know the matrix $D^{-1}$ from \cite[Theorem 9.2]{Spencer_2023}. Moreover, if $C$ is the Cartan matrix for $TL_m$, then $C=D^TD$ (see \cite[Theorem 3.7]{lehrer1996cellular}), so we also obtain the projective characters. 
\end{Remark}

\begin{Example}\label{TL example}
Let $k=\C$. The characters for the cell modules, the simple modules, and the projective modules of $TL_7$ are given in \autoref{TL table} below.  

\begin{table}[H]
\begin{minipage}[b]{0.32\textwidth}
\centering
\begin{tabular}{c|cccc}
$i/j$ & 1 & 3 & 5 & 7 \\\hline
1     & 1 & 2 & 5 & 14\\
3     & 0 & 1 & 4 & 14\\
5     & 0 & 0 & 1 & 6 \\
7     & 0 & 0 & 0 & 1
\end{tabular}
\end{minipage}%
\hfill
\begin{minipage}[b]{0.32\textwidth}
\centering
\begin{tabular}{c|cccc}
$i/j$ & 1 & 3 & 5 & 7 \\\hline
1     & 1 & 1 & 1 & 1 \\
3     & 0 & 1 & 4 & 13\\
5     & 0 & 0 & 1 & 6 \\
7     & 0 & 0 & 0 & 1
\end{tabular}
\end{minipage}%
\quad
\begin{minipage}[b]{0.32\textwidth}
\centering
\begin{tabular}{c|cccc}
$i/j$ & 1 & 3 & 5 & 7 \\\hline
1     & 1 & 2 & 5 & 14 \\
3     & 1 & 3 & 9 & 28\\
5     & 0 & 0 & 1 & 6 \\
7     & 1 & 3 & 9 & 29
\end{tabular}
\end{minipage}%
\caption{From left to right: the character table of cell modules, simple modules, and projective modules.}
\label{TL table}
\end{table}
Note that $i=5$ is critical and so $V_5\cong S_5\cong P_5$. We also have $S_7\cong V_7$ and $S_1\cong P_1$. 

The matrix $L^{-1}$ is obtained by taking transpose of the simple character table and then inverting: $$L^{-1}= \begin{pmatrix}
1 & 0 & 0 & 0\\
-1 & 1 & 0 & 0\\
3 & -4 & 1 & 0\\
-6 & 11 & -6 & 1
\end{pmatrix}$$
with column sums $-3,8,-5,1$. Now if e.g. $V=V_3$, we obtain $l(n)=13^n-5\cdot 4^n+8$ from \autoref{proj length}.
\end{Example}

\begin{Proposition}\label{TL standard ratio}
If $V=S_i$, then $\chi_{\textup{sec}}=\alpha_{m-2,i}$, and so \begin{enumerate}
\item $\lvert l(n)/k(n)-1\rvert \in \mathcal{O}\Big(\big(\frac{(m-i)(m+i+2)}{4m(m-1)}\big)^n\Big),$ and
\item $\lvert l(n)-k(n)\rvert \in \mathcal{O}\big((\alpha_{m-2,i})^n \big).$
\end{enumerate}
\end{Proposition}
\begin{proof}
It is clear that $\chi_{sec}=\alpha_{m-2,i}$ since $\alpha_{j,i}$ is strictly increasing in $j$ when $j\ge i$. The other claims follow from \autoref{proj length}. 
\end{proof}
\begin{Remark}
We can easily write down an analog of \autoref{TL standard ratio} for the $P_i$'s using \autoref{TL irr and proj char}, since $\phi_i(j)$ is again strictly increasing in $j$. 
\end{Remark}

\begin{figure}[H]
\centering
\begin{minipage}{0.450\textwidth}
\centering
\includegraphics[width=0.9\linewidth]{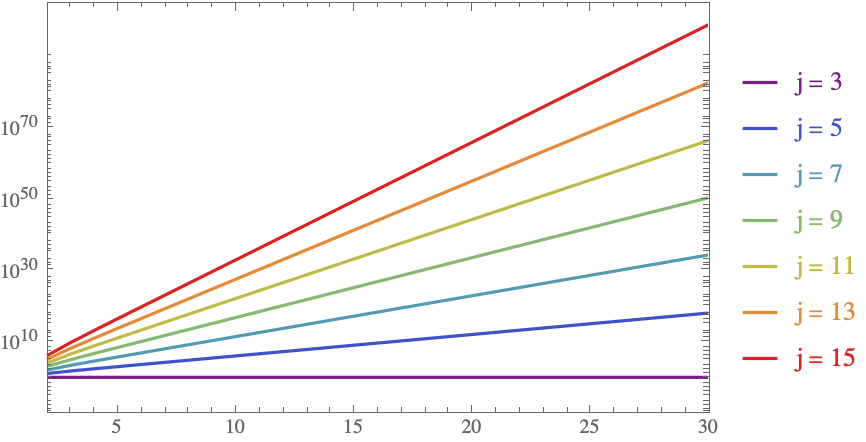}
\end{minipage}
\begin{minipage}{0.45\textwidth}
\centering
\includegraphics[width=0.9\linewidth]{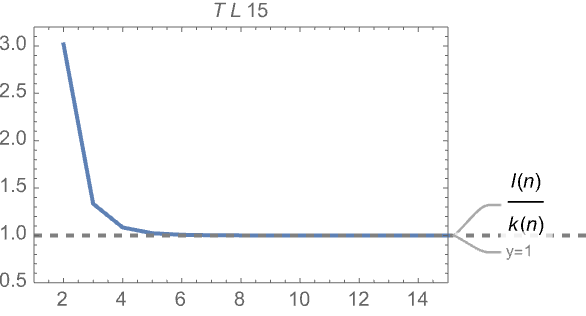} 
\end{minipage}
\caption{Left: the multiplicities $[V^{\otimes n}:V_j]$ for the $\C TL_{15}$-module $V=S_3$ . Right: the ratio $l(n)/k(n)$.}
\label{fig:TL15}
\end{figure}

We illustrate in \autoref{fig:TL15} the case $V=S_3$ in $TL_{15}$, with the left panel showing the multiplicities $[V^{\otimes n}: V_j]$ and the right panel showing the convergence $b(n)/a(n)\to 1$. By \autoref{TL standard ratio}, in this case we have $\lvert l(n)/k(n)-1\rvert \in \mathcal{O}((2/7)^n)$ (which explains the relatively fast rate of convergence), and $\lvert l(n)-k(n)\rvert \in \mathcal{O}(572^n)$. The characters of $V$ are 2002, 572, 165, 48, 14, 4, 1, and 0, with $\chi_{\sec}=572$.

\subsubsection{Motzkin monoid}
For this section we assume $k=\C$. The Motzkin monoid algebra $\C Mo_m$ is isomorphic to the dilute Temperley--Lieb algebra $dTL_m(0)$, and over $\C$ the representation theory of the latter is almost identical to that of the Temperley--Lieb algebra (see e.g.\cite{Bellet_te_2017}). We say $i$ is a critical integer if $i$ is odd, and if $i$ is even, set $i^-=i-2$, $i^+=i+2$. In this notation, the short exact sequences \autoref{eqn: SES ir TL} and \autoref{eqn:SES proj TL} hold verbatim (again with the caveat that $V_{i^-}$ or $V_{i^+}$ is zero if $i^-$ or $i^+$ is out of range). We will use the same strategy to determine $l(n)$ as for the Temperley--Lieb monoid, by inverting the `cell character table'.

\begin{Proposition}
The characters of the cell module $S_i$ is 
\[\chi_{S_i}(j)=\beta_{j,i}:=\begin{cases} \sum_{t=0}^j \frac{i+1}{i+t+1}\binom{j}{i+2t}\binom{i+2t}{t} & j\ge i, \\ 0 & else.\end{cases}\]
\end{Proposition}
\begin{proof}
By an argument identical to that given in the proof of \autoref{TL char}, we have $\chi_{S_i}(j)=\dim S_i$ where $S_i$ is considered as a $\C Mo_j$-module. The result now follows by computing the dimensions of the cell modules, see \cite[Proposition 4F.4]{KhSiTu-monoidal-cryptography}.
\end{proof}
\begin{Remark} \label{Mo char remark}
Again, these character values are independent of $m$, so the cell character table for each $Mo_m$ is a truncation of an infinite uni upper triangular matrix $S$ with rows and columns indexed by $\N$, $S_{ij}=\beta_{j,i}$. In other words, the $i$th row of $S$ corresponds (up to an index shift) to the $i$th convolution of the Motzkin numbers.  
\end{Remark}
\begin{Lemma}\label{inverse Mo}
With $S$ as in \autoref{Mo char remark}, the inverse of $S$ is $$S^{-1}_{ij}=\begin{cases}
[x^{j-i}](1+x+x^2)^{-(i+1)} = (-1)^{j-i}\sum_{r=0}^{[\frac{j-i}{2}]}(-1)^r\binom{i+r}{r}\binom{j-r}{j-i-2r} & i\le j,\\
0 & i>j.
\end{cases}$$
\end{Lemma}
\begin{proof}
We know that the $i$th row of $S$ is the $i$th convolution of the Motzkin numbers. Equivalently, $S$ is the Riordan array $(M(x),xM(x))$, where $M(x)=1+x+2x^2+4x^3+\dots$ is the usual generating function for the Motzkin numbers satisfying $M(x)=1+xM(x)+x^2(M(x))^2$. The composition inverse of $xM(x)$ is $x/(1+x+x^2)$, and so $S^{-1}$ is exactly the inverse Riordan array $(1/({1+x+x^2}),x(1+x+x^2)).$ The claimed formula now follows by computation. 
\end{proof}

\begin{Proposition}
Let $V$ be a $kMo_m$-module with character $\chi$. With $S^{-1}$ as above, we have
$$[V^{\otimes n}: V_i]=\begin{cases}
\sum_{j} (\chi(j))^n S^{-1}_{ji} &  i \text{ odd or $i=0$,}\\
\sum_{j} (\chi(j))^n S^{-1}_{ji} + \sum_{j} (\chi(j))^n S^{-1}_{ji^{-}} & \text{ else}.\end{cases}$$
Summing all the terms give a formula for $l(n)$.
\end{Proposition}
\begin{proof}
As in the proof of \autoref{TL length}, we can write $\chi$ in the basis of the characters of cell modules, so $$\chi = \sum_i c_i\chi_{S_i}, c_i=\sum_{j}(\chi(j))S_{ji}^{-1},$$ and the result now follows from the short exact sequence \autoref{eqn: SES ir TL}.
\end{proof}

We can also compute the simple and projective character tables for $Mo_m$.  Again let $\chi_i$ (resp. $\phi_i$) denote the character of $V_i$ (resp. $P_i$).
\begin{Proposition}
If $i$ is odd, then $\chi_i=\chi_{S_i}$. If $i=0$, then $\chi_i$ is the trivial character. If $i=2l, l>0$, then 
we have 
\begin{equation}\label{motz irr}
\chi_{2l}(j)=
\sum_{\substack{0\le t\le j-2l\\ t\equiv j\pmod 2}}
\frac{4l}{j-t+2l}\,
\binom{j}{t}\,
\binom{j-t-1}{(j-t)/2+l-1}.
\end{equation}
If $i$ is odd or $i=0$, then $\phi_i=\chi_{S_i}$. Otherwise, $\phi_i(j)=\beta_{j,i}+\beta_{j,i-2}.$ 
\end{Proposition}
\begin{proof}
The statements about $\phi_i$ are clear from \autoref{eqn:SES proj TL}. From \autoref{eqn: SES ir TL}, we also have $$\chi_{{m'-2i}}(j)=\sum_{s=0}^m (-1)^{s}\beta_{j,m'-2i+2s},$$ where $m'=m$ if $m$ is odd, and write $m'=m-1$ if $m$ is even. Writing $2l=m'-2i$, it can be verified (by induction or otherwise) that this formula is equivalent to \autoref{motz irr}.  
\end{proof}

\begin{Remark}
The formula in \autoref{motz irr} correspond to columns of the triangle $T(n,k)$ given in \cite[A379838]{Oeis}, which counts the total number of humps with height $k$ in all Motzkin paths of order $n$.
From this formula we immediately get that the character table of $Mo_m$ has only positive entries the upper triangular part, and in fact $\chi_i(j)>\chi_{i}(l)\ge 1$ for $j>l\ge i$.
\end{Remark}

\begin{Proposition} 
If $\chi$ is the character of $V$, where either $V=S_i$, $V=P_i$, or $i>0$ and $V=V_i$, then we have $\chi_{\textup{sec}}=\chi(m-1).$ 
\end{Proposition}
\begin{proof}
In all cases $\chi(j)$ is strictly increasing when $j\ge i>0$. 
\end{proof}

\begin{figure}[H]
\centering
\begin{minipage}{0.450\textwidth}
\centering
\includegraphics[width=0.9\linewidth]{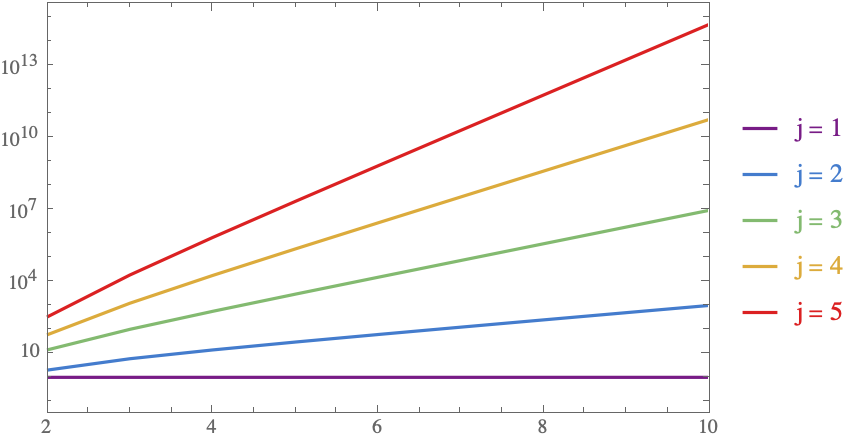}
\end{minipage}
\begin{minipage}{0.45\textwidth}
\centering
\includegraphics[width=0.9\linewidth]{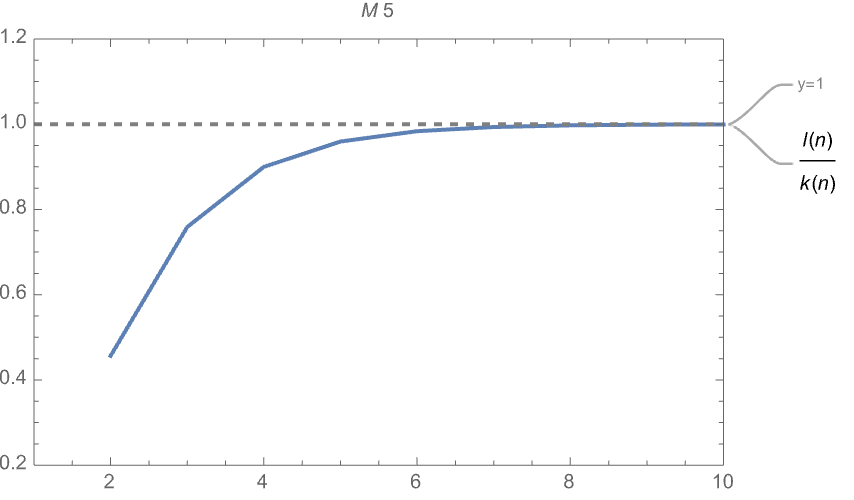} 
\end{minipage}
\caption{Left: the multiplicities $[V^{\otimes n}:V_j]$ for the $\C Mo_5$-module $V=S_1$. Right: the ratio $l(n)/k(n)$.}
\label{fig:M5}
\end{figure}

\begin{Example}\label{Mo example}
The characters for the cell modules, the simple modules, and the projective modules (over $\C$) of $Mo_5$ are given in \autoref{M table} below.  

\begin{table}[H]
\begin{minipage}[b]{0.32\textwidth}
\centering
\begin{tabular}{c|cccccc}
$i/j$ & 0 & 1 & 2 & 3 & 4& 5 \\\hline
0     & 1 & 1 & 2 & 4 & 9 & 21\\
1     & 0 & 1 & 2 & 5 & 12 & 30\\
2 & 0 & 0 & 1 & 3 & 9 & 25\\
3     & 0 & 0 & 0 & 1 & 4 & 14 \\
4     & 0 & 0 & 0 & 0 &1 & 5 \\
5 & 0 & 0 & 0 & 0 & 0  &1
\end{tabular}
\end{minipage}%
\hfill
\begin{minipage}[b]{0.32\textwidth}
\centering
\begin{tabular}{c|cccccc}
$i/j$ & 0 & 1 & 2 & 3 & 4 & 5 \\\hline
0 & 1 & 1 & 1 & 1 & 1 & 1 \\
1 & 0 & 1 & 2 & 5 & 12 & 30 \\
2 & 0 & 0 & 1 & 3 & 8 & 20 \\
3 & 0 & 0 & 0 & 1 & 4 & 14 \\
4 & 0 & 0 & 0 & 0 & 1 & 5 \\
5 & 0 & 0 & 0 & 0 & 0 & 1
\end{tabular}

\end{minipage}%
\quad
\begin{minipage}[b]{0.32\textwidth}
\centering
\begin{tabular}{c|cccccc}
$i/j$ & 0 & 1 & 2 & 3 & 4& 5 \\\hline
0     & 1 & 1 & 2 & 4 & 9 & 21\\
1     & 0 & 1 & 2 & 5 & 12 & 30\\
2 & 1 & 1 & 3 & 7 & 18 & 47\\
3     & 0 & 0 & 0 & 1 & 4 & 14 \\
4     & 0 & 0 & 0 & 0 &10 & 30 \\
5 & 0 & 0 & 0 & 0 & 0  &1
\end{tabular}
\end{minipage}%
\caption{From left to right: the character table of cell modules, simple modules, and projective modules.}
\label{M table}
\end{table}
Note that all characters are strictly increasing in $j$ to the right of the main diagonal. 
We obtain $L^{-1}$ by taking the inverse transpose of the simple character table: $$L^{-1}=
\begin{psmallmatrix}
1 & 0 & 0 & 0 & 0 & 0 \\
-1 & 1 & 0 & 0 & 0 & 0 \\
0 & -2 & 1 & 0 & 0 & 0 \\
1 & 1 & -3 & 1 & 0 & 0 \\
-1 & 2 & 3 & -4 & 1 & 0 \\
0 & -4 & 2 & 6 & -5 & 1
\end{psmallmatrix}
$$ with column sums $0,-2,3,3,-4,1$. So for e.g.$V=S_1$, we have $l(n)=-2+3\cdot 2^n +3\cdot 5^n  -4\cdot 12^n + 30^n$.

We illustrate the $V=S_1$ case below in \autoref{fig:M5}, with the left panel showing the growth rate of the composition multiplicities and the right panel showing the convergence of the ratio $l(n)/k(n)$ to 1. In this case $\chi_{\mathrm{sec}}=12$, so we have $\lvert l(n)/k(n)-1\rvert \in \mathcal{O}((2/5)^n),$ and $\lvert l(n)-k(n)\rvert \in \mathcal{O}(12^n).$ The characters of $V$ are entries in the second row in \autoref{M table}.
\end{Example}

\subsection{Length for full and special linear monoids}
We now consider an example where $k(n)$ is nontrivial. Denote by $\M(2,q)$ the full linear monoid of $2 \times 2$ matrices with entries in the finite field $\mathbb{F}_q$, where $q=p^r$ is a prime power, and similarly denote by $\SLM(2,q)$ the special linear monoid of $2\times 2$ matrices in $\mathbb{F}_q$ with determinant 0 or 1. The group of units for $\M(2,q)$ and $\SLM(2,q)$ are the general and special linear groups $\GL(2,q)$ and $\SL(2,q)$ respectively. 

\begin{Proposition}\label{matrix monoid}
Let $k$ be a field of characteristic $p>0$, and let $q=p^r$. Let $G$ denote either $\textup{SL}(2,q)$ or $\textup{GL}(2,q)$ and let $M$ denote either $\textup{M}(2,q)$ or $\textup{SLM}(2,q)$. Suppose that $V$ is either a faithful $kG$-module, or a $kM$-module on which no element apart from 1 acts as identity, then

$$l(n)\sim k(n)=\frac{2r-1+(2p-1)^r-2^r}{q^2-1} \cdot (\dim V)^n. $$ 
\end{Proposition}
\begin{proof}
Under the given assumptions it suffices to show that for a faithful $kG$-module the formula for $k(n)$ is as stated. For $G=\SL(2,q)$, we apply \autoref{proj length} using the projective Brauer character table of $\textup{SL}(2,q)$ in defining characteristic worked out in \cite[\S 4]{Sri}. For $\textup{GL}(2,q)$ there are $(q-1)$ times more characters coming from tensoring with some power of the determinant representation, but the order of the group also increases by a factor of $(q-1)$, so $k(n)$ remains unchanged. 
\end{proof}
\begin{Remark}
The formula for $k(n)$ in the case when $\M(2,q)$ is semisimple is worked out in \cite[Proposition 4B.1]{he2025growthproblemsrepresentationsfinite}.
Moreover, the formula for $k(n)$ given above in \autoref{matrix monoid} is aperiodic, but the corresponding formula for $a(n)$ is periodic of period 2, see \cite[Proposition 4B.3]{he2025growthproblemsrepresentationsfinite}.
\end{Remark}

\section{Growth rate of the number of summands of tensor powers}

\subsection{General theory}
In \cite[Conjecture 1B.1]{he2025growthproblemsrepresentationsfinite}, it was conjectured that if $V$ is a $kM$-module on which no element apart from $1\in M$ acts as identity, then 
\begin{equation}\label{eqn: bn conj}
b(n) \sim a(n)=\frac{1}{\lvert G\rvert}\sum_{\substack{1\le t\le N \\ g_t \in Z_V(G)}} S_{t}\big(\omega_V(g_t)\big)^n\cdot (\dim V)^n,    
\end{equation}
where $1\le t\le N$ index the conjugacy classes of the group of units $G$, and $S_t$ is the sum over the conjugate of the $t$-th column of the simple (Brauer) character table of $G$. This is the analog of \autoref{eqn:no identity length formula} for $b(n)$, and states that $a(n)$ for a $kM$-module $V$ is the same as $a(n)$ for the $kG$-module $\Res(V)$ obtained by restriction. (The formula for $a(n)$ of a faithful $kG$-module is well understood from \cite[Theorem 1]{he2024growthproblemsrepresentationsfinite}.) Given a $kG$-module $V$, we also obtain a $kM$-module $\Ind(V)$ by letting the nonunits act as 0. The conjecture above is motivated by the following conjecture, whose truth would imply it:
\begin{Conjecture}\label{conj}
Let $V$ be a $kM$-module on which no element apart from $1\in M$ acts as identity, then there is some $n\in N$ such that $V^{\otimes n}$ contains a direct summand of the form $\Ind(W)$ for a $kG$-module $W$.
\end{Conjecture}
To see that this implies \autoref{eqn: bn conj}, first observe that if $V$ is a $kM$-module and $W$ is a $kG$-module, then 
\begin{equation} \label{Ind res} V\otimes \Ind(W)=\Ind(\Res(V)\otimes W).
\end{equation}
If the conjecture is true, then there is some $n_0\in N$ such that for $n\ge n_0$, $V^{\otimes n}$ contains a summand of the form $\Ind(\Res(V^{\otimes (n-n_0)})\otimes W)$. By assumption, $\Res(V)$ is a faithful $kG$-module, and so by \cite[Theorem 2]{bryant1972tensor} for $n$ large enough $\Res(V^{\otimes (n-n_0)})$ has a projective summand $P$, so $V^{\otimes n}$ has a summand of the form $\Ind(Q)$ where $Q$ is a projective $kG$-module. Define the (potentially countably infinite) fusion graph (for summands) $\Gamma_b$ associated with $V$ by the directed and weighted graph whose vertices are indecomposable $kM$-modules that occur as a summand in some tensor power $V^{\otimes n}$, such that we draw an edge of weight $m$ from $W$ to $X$ if $V\otimes W$ has $m$ copies of $X$ in the Krull--Schmidt decomposition. Noting that modules of the form $\Ind(P)$ ($P$ projective) form a tensor ideal, the above discussion shows that there is a path from every vertex in $\Gamma_b$ to the strongly-connected component $\Gamma^P_b$ of $\Gamma_b$ induced by vertices (finitely many) which are inductions of projective indecomposable $kG$-modules; moreover, $\Gamma_b^G$ has no outgoing edges. This now implies that asymptotically, $l(n)$ is controlled only by $\Gamma_b^G$ (for the technical details we refer to \cite[Theorem 3A.2] {he2025growthproblemsrepresentationsfinite} and \cite[Proposition 4.22]{lacabanne2024asymptotics}). Indeed, \autoref{conj} would also imply bounds on the rate of convergence and variance, as well as an asymptotic formula for the number of times an indecomposable $X$ appears as a summand in $V^{\otimes n}$, in analog with the formulas in \autoref{proj length} (everywhere $\chi_{\sec}$ should be replaced by the \textit{second largest eigenvalue} of the adjacency matrix of $\Gamma_b$). 

\begin{Definition}\label{def}
We say $kM$ satisfies the group-injective condition if there is some projective indecomposable $kG$-module $P$ such that $\Ind(P)$ is injective as a $kM$-module. 
\end{Definition}

\begin{Theorem}\label{T:Main}
The group-injective condition implies \autoref{eqn: bn conj} and \autoref{conj}.
\end{Theorem}

\begin{proof}
By \cite[Theorem 3A.2]{he2025growthproblemsrepresentationsfinite}, we know that \autoref{eqn: bn conj} and indeed \autoref{conj} holds when the group-injective condition is satisfied. 
\end{proof}

In the group-injective case we have an analog of \autoref{n0 bound}:

\begin{Proposition} \label{m0 bound}
Let $m_0$ denote the smallest natural number such that $V^{\otimes m_0}$ has a summand of the form $\Ind(P)$, where $P$ is a projective indecomposable $kG$-module. If $M$ satisfy the group-injective condition, and $V$ is such that no element apart from $1$ acts as identity on it, then
$$m_0\le L+\frac{\lvert G\rvert}{\lvert Z_{\Res(V)}(G)\rvert}+\lvert Z_{\Res(V)}(G)\rvert-3,$$
where $L$ is the number of $\mathcal{L}$-classes of $M$.
\end{Proposition}
\begin{proof}
Write $g:=\lvert G\rvert$ and $m:=\lvert Z_{\Res(V)}(G)\rvert.$ By the proof of \autoref{n0 bound}, we know that $V^{\otimes(L-1)}$ contains a submodule of the form $\Ind(W)$, for $W$ a $kG$-module. By \cite[Theorem 2]{bryant1972tensor}, the regular module $kG$ is a summand of $$(\Res(V))^{\otimes (g/m-1)} \oplus (\Res(V))^{\otimes (g/m)}\oplus \dots \oplus (\Res(V))^{\otimes (g/m+m-2)}.$$
Since $W\otimes kG=kG^{\oplus (\dim W)}$, by \autoref{Ind res} we see that $V^{\otimes (L+g/m-2)}\oplus \dots \oplus V^{\otimes (L+g/m+m-3)}$ contains a submodule of the form $\Ind(kG)$, but by the group-injective condition, we now have that one of the modules $V^{\otimes (L+g/m-1)},\dots, V^{\otimes (L+g/m+m-3)}$ contains a summand of the form $\Ind(P)$ where $P$ is projective, so we are done.
\end{proof}

We will now establish that various families of monoids satisfy the conjectured equation \autoref{conj}, by showing that they satisfy the group-injective condition or otherwise. Among other things, we will prove:
\begin{Theorem}\label{diagram bn}
Let $k=\C$, and let $V$ be a $\C M$-module on which no element apart from $1\in M$ acts as identity. Let $m\ge 1$.
\begin{enumerate}
\item If $M=pRo_m$, $TL_m$ or $Mo_m$, then \begin{equation} \label{eqn: planar asym} a(n)\sim (\dim V)^n.\end{equation}
\item If $M=Ro_m, Br_m, RoBr_m$, or $Pa_m$, then  
\begin{equation} \label{eqn: sym asym} a(n)\sim \sum_{z=0}^{\floor * {m/2}} \frac{1}{(m-2z)!z!2^z}\cdot (\dim V)^n.\end{equation}
\end{enumerate}
\end{Theorem}
\begin{proof}
The $M=Ro_m$ case follows from \autoref{plannarrook formula}. The $M=Ro_m$ case is known from \cite[Proposition 4A.1]{he2025growthproblemsrepresentationsfinite}; see also \autoref{rook cor}. For $M=Br_m, RoBr_m$ or $Pa_m$ see \autoref{sym monoids cor}. For $M=TL_m$ or $Mo_m$, see \autoref{TL Mo an}. 
\end{proof}

\subsection{Monoids satisfying the group-injective condition}

We will use \autoref{T:Main} silently throughout this section.

\subsubsection{Inverse monoids}
\begin{Proposition}
Let $M$ be an inverse monoid, then over $k$ of any characteristic the group-injective condition holds.    
\end{Proposition}
\begin{proof}
Recall (e.g.from \cite[Corollary 9.4]{steinberg2016representation}) that the monoid algebra $kM$ for an inverse monoid $M$ is isomorphic to the direct sum of certain matrix algebras, one summand of which being isomorphic to $kG$. It follows that the induction of any projective $kG$-module remains projective as a $kM$-module. Finally, for an inverse monoid the projective and injective modules coincide (see \cite[\S 15.2]{steinberg2016representation}).
\end{proof}

This applies to, for example, the rook monoid $Ro_m$ (also called the symmetric inverse monoid), which can be viewed as the submonoid of the partition monoid where all blocks have size at most $2$, such that every block of size $2$ is \textit{propagating}, i.e. consists of one top vertex and one bottom vertex. The group of units of $Ro_m$ is the symmetric group $S_m$, so we obtain:
\begin{Proposition}\label{rook cor}
Let $V$ be a $kRo_m$-module with no elements apart from $1$ acting as identity, then 
$$a(n)=\frac{(\dim V)^n}{m!} \sum_{S \in \textup{Irr}_k(S_m)} \dim S,$$
where the sum is over the isomorphism classes of simple $kS_m$-modules.
\end{Proposition}

\begin{proof}
The group of units of $Ro_m$ is the symmetric group $S_m$, and the result follows.
\end{proof}

In nonmodular characteristics ($p\nmid n!$), the sum of dimensions of simple $kS_m$-modules is known, see, for example, \cite[Example 2.3]{coulembier2023asymptotic} or \cite[A000085]{Oeis}.

\begin{figure}[H]
\centering
\begin{minipage}{0.40\textwidth}
\centering
\includegraphics[width=1\linewidth]{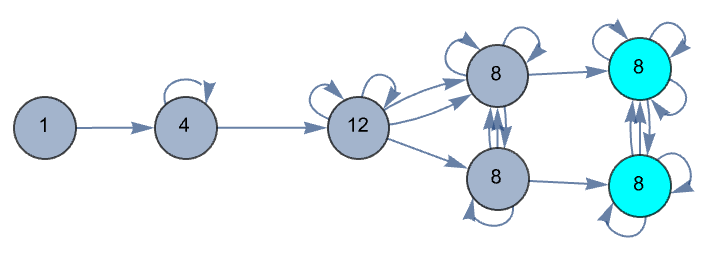}
\end{minipage}
\begin{minipage}{0.4\textwidth}
\centering

\[
\begin{psmallmatrix}
0 & 0 & 0 & 0 & 0 & 0 & 0 \\
1 & 1 & 0 & 0 & 0 & 0 & 0 \\
0 & 1 & 2 & 0 & 0 & 0 & 0 \\
0 & 0 & 2 & 2 & 2 & 0 & 0 \\
0 & 0 & 1 & 1 & 1 & 0 & 0 \\
0 & 0 & 0 & 1 & 0 & 3 & 2 \\
0 & 0 & 0 & 0 & 1 & 1 & 2
\end{psmallmatrix}
\]

\end{minipage}
\caption{The fusion graph and its adjacency matrix for the four-dimensional simple $\mathbb{F}_2Ro_4$-module. The component $\Gamma_b^G$ is colored in cyan.}
\label{fig:rookmonoid}
\end{figure}

\begin{figure}[H]
\centering
\includegraphics[width=0.4\linewidth]{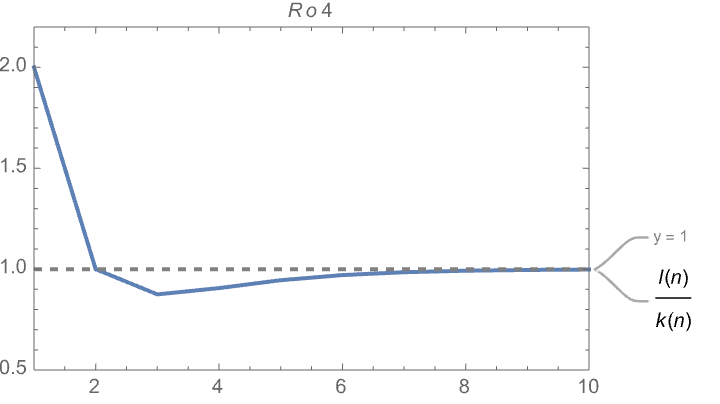}
\caption{The ratio $b(n)/a(n)$ for the four-dimensional simple $\mathbb{F}_2Ro_4$-module.}
\label{fig:I4_conv}
\end{figure}

We illustrate in \autoref{fig:rookmonoid} the fusion graph (and its adjacency matrix) in the case when $V$ is the four-dimensional simple representation of $Ro_4$ in characteristic 2. We have drawn the labeled edges as multi-edges and labeled the nodes by their dimensions.  The strongly-connected component $\Gamma_b^G$, colored in cyan, consist of (induction of) the two projective indecomposables for $S_4$, both of dimension 8. It takes four steps for the trivial module to reach $\Gamma_b^G$, and so  $m_0=4$. In characteristic 2 the group $S_4$ has two simple representations whose degrees sum up to 3, and so $a(n)=1/8\cdot 4^n$. The ratio $b(n)/a(n)$ is illustrated in \autoref{fig:I4_conv}.

\subsubsection{Symmetric monoids}
We now show that the Brauer monoid $Br_m$, rook-Brauer monoid $RoBr_m$ and partition monoid $Pa_m$ satisfy the group-injective condition for all $m$, when $k=\C$. Indeed, we will show that the module $\Ind(V_{sgn})$, where $V_{sgn}$ is the sign representation of $S_m$, is injective in $\C M$ where $M$ is $Br_m, RoBr_m$ or $Pa_m$. We will show this by exhibiting an idempotent $e$ such that $\C Me \cong \Ind(V_{sgn})$. 

Recall that $Br_m$ is the submonoid of $Pa_m$ consisting of all diagrams where each block has size exactly 2. An element in $Br_m$ is invertible if and only if each block consists of one top vertex and one bottom vertex (that is, there are no `cups' or `caps'). It is clear that the invertible elements may be identified with $S_m$.   
\begin{Proposition}\label{Brauer}
Let $e=\frac{1}{m!}\sum_{\sigma \in S_m} (sgn(\sigma))\sigma$ be the usual idempotent for $V_{sgn}$ in $S_m$, then $\C Br_m e \cong \Ind(V_{sgn})$ and so $\C Br_m$ satisfies the group-injective condition.
\end{Proposition}
\begin{proof} To see that $\C Br_m e \cong \Ind(V_{sgn})$, it suffices to show that if $x \in Br_m \setminus S_m$, then $xe=0$. Let $x\in Br_m \setminus S_m$, then it must have a `cap' joining the bottom vertices with labels $i$ and $j$. It follows that if we represent the transposition $(i,j) \in S_m$ by the obvious diagram in $Br_m$ with a single $(i,j)$-crossing, then $x (i,j)=x$. Thus we have $xe=x(i,j)e=-xe$, so $xe=0$. This shows that $\Ind(V_{sgn})$ is projective. Let $\iota$ denote the anti-isomorphism of $Br_m$ which reflects the diagrams across a horizontal axis, then by equipping $D(V)=\textup{Hom}(V,k)$ with the action $x\cdot f(v)=f(\iota(x) v)$, we get a simple-preserving duality $V\mapsto D(V)$ sending projective modules to injective modules (see e.g. \cite[\S 2.3]{Bellet_te_2017} for the argument), so $\Ind(V_{sgn})\cong D(\Ind(V_{sgn}))$ is also injective.
\end{proof}

Now we consider $RoBr_m$ and $Pa_m$. We have the same embedding of $S_m$ into these monoids, though now a nonunit could have isolated vertices (`dots') in addition to cups and caps. In $RoBr_m$ the blocks are required to have size at most 2, whereas in $Pa_m$ the blocks can have any size. For $1\le i \le m$, denote by $P_i$ the diagram in $RoBr_m$ or $Pa_m$ consisting of two dots in the $i$th position and vertical lines elsewhere. If $\sigma \in S_m$, it is clear that $p_i\sigma = \sigma p_{\sigma(i)}$.

\begin{Proposition}\label{Robr idem}
Let $M=RoBr_m$ or $Pa_m$, and let $$e=\frac{1}{m!}\big(\sum_{\sigma \in S_m} \big(sgn(\sigma))\sigma+\sum_{i=1}^m \sum_{\sigma \in S_m} \Big((-sgn(\sigma))\sigma p_i\Big).$$ Then $\C Me\cong \Ind(V_{sgn})$, and $\C M$ satisfies the group-injective condition. 
\end{Proposition}
\begin{proof}
It is clear that if $\sigma \in S_m$, then $\sigma e=sgn(\sigma)e$, so it suffices to see that $xe=0$ for $x\notin S_m$. Suppose $x\notin S_m$. If $m$ has two connected points $i,j$ on the bottom, then $x(i,j)=x$ and the same argument as in the proof of \autoref{Brauer} shows that $x \sum_{\sigma\in S_m}(-sgn(\sigma))\sigma p_i=0$ for each $i$, so $xe=0$. Suppose $x$ has no connected points on the bottom, then it must have an isolated point, say at the $k$th position, and so $xp_k=x$. As $p_i\sigma = \sigma p_{\sigma(i)} ,$ we can also write $e=\frac{1}{m!}\big(\sum_{\sigma \in S_m} (sgn(\sigma))\sigma+\sum_{i=1}^n \sum_{\sigma \in S_m} (-sgn(\sigma))p_i\sigma \big)$, whence $$x\sum_{\sigma \in S_m}(sgn(\sigma))\sigma + x\sum_{\sigma \in S_m}(-sgn(\sigma)p_k\sigma)=0.$$ On the other hand, if $i\neq k$ then $xp_i$ has at least two isolated points, $i$ and $k$, and $xp_i(i,k)=xp_i$. So $$x\sum_{\sigma \in S_m}(-sgn(\sigma))p_i\sigma=x\sum_{\sigma \in S_m}(-sgn(\sigma))p_i (i,k)\sigma = -x\sum_{\sigma \in S_m}(-sgn(\sigma))p_i\sigma$$ so the left hand side is $0$ for each $i\neq k$. Thus, we have $\C Me\cong \Ind(V_{sgn})$, and by the same duality argument as in the proof of \autoref{Brauer} $\Ind(V_{sgn})$ is also injective.
\end{proof}

\begin{Proposition}\label{sym monoids cor}
Let $V$ be a $\C M$-module on which no element apart from $1\in M$ acts as identity, where $M=Br_m,RoBr_m$ or $Pa_m$, then we have 
$$a(n)=\sum_{z=0}^{\floor * {m/2}} \frac{1}{(m-2z)!z!2^z}\cdot (\dim V)^n.$$
\end{Proposition}
\begin{proof}
By \autoref{Brauer} and \autoref{Robr idem}, the monoids involved all satisfy the group-injective condition over $\C$. The result now follows from \autoref{eqn: bn conj} and the computation that one has $\sum_{S \in \textup{Irr}_k(S_m)} \dim S=\sum_{z=0}^{\floor * {m/2}} 1/{(m-2z)!z!2^z}$ (see \cite[Example 2.3]{coulembier2023asymptotic}.)
\end{proof}

\begin{Proposition}\label{idem multiplicity}Let $V$ be a $\C M$-module with character $\chi$. 
\begin{enumerate}
\item If $M=Br_m$, we have $$[V^{\otimes n}:\Ind(V_{\text{sgn}})]=[V^{\otimes n}:\Ind(V_{sgn})]_b=\frac{1}{n!}\sum_{\sigma \in S_m} sgn(\sigma)(\chi(\sigma))^n,$$ depending only on the restriction of $V$ to $S_m$. If $\Res(V)$ is faithful, then $[V^{\otimes n}: V_{sgn}]\sim \frac{1}{m!}(\dim V)^n$. 
\item If $M=RoBr_m$ or $Pa_m$, we have $$[V^{\otimes n}: \Ind(V_{sgn})]=[V^{\otimes n}:\Ind(V_{sgn})]_b= \frac{1}{m!}\Big( \sum_{\sigma \in S_m} sgn(\sigma) (\chi(\sigma))^n + \sum_{i=1}^m \sum_{\sigma \in S_m} -sgn(\sigma) (\chi(\sigma p_i))^n\Big).$$
\end{enumerate}
\end{Proposition}
\begin{proof}
That $[V^{\otimes n}: \Ind(V_{sgn})]$ is given by the claimed formula follows from Propositions \ref{Brauer}, \ref{Robr idem}, and the well-known fact (see e.g. \cite[Proposition 7.19]{steinberg2016representation}) that if $e=\sum_{m\in M}c_mm$ is the primitive idempotent corresponding to a simple module $S$, then $$[V:S]=\sum_{m \in M}{c_m}\chi(m).$$ Finally, note that since $\Ind(V_{sgn})$ is both projective and injective, it occurs as a composition factor if and only if it occurs as a summand. 
\end{proof}

\begin{figure}[H]
\centering
\begin{minipage}{0.40\textwidth}
\centering
\includegraphics[width=0.6\textwidth]{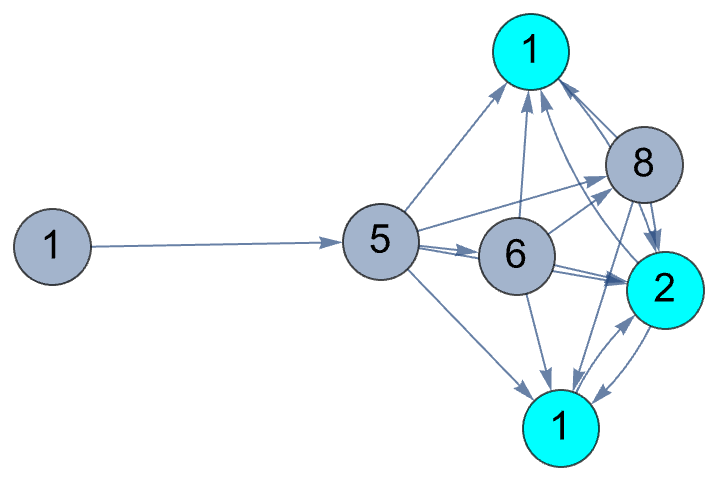} 
\end{minipage}
\begin{minipage}{0.4\textwidth}
\centering

\[
\begin{psmallmatrix}
0 & 0 & 0 & 0 & 0 & 0 & 0 \\
1 & 1 & 0 & 0 & 0 & 0 & 0 \\
0 & 1 & 3 & 0 & 1 & 2 & 3 \\
0 & 3 & 0 & 3 & 1 & 4 & 5 \\
0 & 1 & 1 & 1 & 4 & 5 & 8 \\
0 & 1 & 0 & 0 & 0 & 1 & 0 \\
0 & 1 & 0 & 0 & 0 & 1 & 2
\end{psmallmatrix}
\] 

\end{minipage}
\caption{The fusion graph and adjacency matrix for the five-dimensional projective $\C Pa_3$-module. Nodes are labeled by dimension, and the component $\Gamma_b^G$ is colored in cyan.}
\label{fig:Pa3}
\end{figure}

We plot in \autoref{fig:Pa3} the fusion graph for the five-dimensional projective $\C Pa_3$-module. Multi-edges and self-loops are suppressed in the graph. Note that since we are in characteristic zero, the projective $\C S_3$-modules are precisely the simple ones. 

\subsubsection{Temperley--Lieb and Motzkin monoid}

Over a field $k$, $TL_m$ or $Mo_m$ satisfies the group-injective condition for some, but not all values of $m$. To state the result recall the notation from \autoref{SS:TL}, and we say $a$ is $(p,l)$ ancestorless if $a=[a_t,0,\dots,0]_{p,l}$. (The nomenclature is a gender neutral variation of the one used in \cite{BuLiSe-tl-char-p,TuWe-quiver-tilting}.) For example, for $p=\infty$ being $(p,l)$ ancestorless means divisible by $l$.

\begin{Proposition}\label{P:GI}
Over an arbitrary field $k$ of characteristic $p$, $TL_m$ satisfies the 
group-injective condition if and only if $m+1$ is $(p,3)$ ancestorless, and $Mo_m$ satisfies the group-injective condition if and only if $m+1$ is $(p,2)$ ancestorless.
\end{Proposition}

\begin{proof}
We now assume familiarity with $(p,l)$ Jones--Wenzl idempotents 
and their relation to tilting modules, cf. \cite{MaSp-pl-jones-wenzl,SuTuWeZh-mixed-tilting}.

For $TL_m$ we use the equivalence of the Temperley--Lieb category with circle parameter $1$ with tilting modules of quantum SL2, see \cite{RuTeWe-sl2}, which happens for $l=1$. The group-injective condition then corresponds to the $(p,l)$ Jones--Wenzl idempotent for $m$ strands to be given by the same formula as the (usual) Jones--Wenzl idempotent for $m$. This happens if and only if when the stated condition is satisfied.

For $Mo_m$ exactly the same argument works but for $l=2$ which is when the Motzkin category is equivalent to tilting modules for quantum SL2, by \cite{BeHa-motzkin}.

(The above references are only in characteristic zero. Making them independent of the field can be done using standard arguments about tilting modules as, for example, in \cite[Proposition 2.3]{AnStTu-semisimple-tilting}; for background on tilting modules see e.g. \cite{Do-tilting-alg-groups,Do-q-schur,Ri-good-filtrations,AnStTu-cellular-tilting}.)
\end{proof}

\subsubsection{Full transformation monoid}

Recall the full transformation monoid $T_m$, which can be defined as the monoid $End(\{1,\dots,m\})$.

\begin{Proposition}
Over $\C$, the monoid $T_m$ satisfies the group-injective condition. Consequently, if $V$ has no element apart from $1$ acting as identity, then $a(n)$ is as in \autoref{eqn: sym asym}.
\end{Proposition}

\begin{proof}
This was already done in \cite{he2025growthproblemsrepresentationsfinite}, and is listed for completeness.
\end{proof}



\subsubsection{Partial transformation monoid and its submonoids}\label{ptm} Let $PT_m$ denote the monoid of all partial functions on $m$ points, which contains as submonoids the monoid $PO_m$ of order preserving partial functions on $m$ points, the partial Catalan monoid $PC_m$, and the monoid $PF_m$ of partial order-decreasing functions. Recall that the (ordinary) quiver of an algebra $A$ is the directed graph whose vertices are the simple $A$-modules, with $\dim \Ext^1(X,Y)$ arrows drawn from $X$ and $Y$. Thus, a simple $A$-module $X$ is injective if and only if there are no arrows ending in $X$ in the quiver of $A$. Over $\C$, the group-injective condition is therefore equivalent to there being no arrow ending in some $\Ind(V_i)$ in the quiver of $\C M$, where $V_i$ is some simple $\C G$-module. 

\begin{Proposition}
Over $\C$, the monoids $PT_m$, $PO_m$, $PC_m$ and $PF_m$ satisfy the group-injective condition. Consequently, if $V$ has no element apart from $1$ acting as identity, then $a(n)$ is as in \autoref{eqn: sym asym} if $M=PT_m$, and $a(n)$ is as in \autoref{eqn: planar asym} if $M=PO_m, PC_m$ or $PF_m$.
\end{Proposition}
\begin{proof} By \cite[Theorem 3.8]{STEIN2016549}, in the quiver of $PT_m$ there are no arrows ending in any of the inductions of $kG$-modules, so for any simple $kG$-module $S$ we have that $\Ind(S)$ is injective. Thus $PT_m$ satisfies the group-injective condition, and the statement about $a(n)$ follows since $G=S_m$ in this case. The proof for $M=PO_m$, $PC_m$ and $PF_m$ are all similar, following from the characterisations of their quivers in \cite[Proposition 5.2; Proposition 5.8; Proposition 5.5]{STEIN2016549}. 
\end{proof}
\subsubsection{Catalan monoid}
Recall that the Catalan monoid $C_m$, $m\ge 1$, is the monoid of all order-preserving and nondecreasing
maps on $m$ points. It is $\mathcal{J}$-trivial and hence $G=1$.
\begin{Proposition}
The Catalan monoid $C_m$ satisfies the group-injective condition over any algebraically closed field $k$.     
\end{Proposition}
\begin{proof}
It is known (see \cite[\S 3.7.3]{denton2011representationtheoryfinitejtrivial} or \cite[Theorem 17.24]{steinberg2016representation}) that the quiver of $kC_m$ is isomorphic to the Hasse diagram of a partial order on the collection $P_m$ of subsets of $\{1,\dots, m\}$ containing $m$, such that $S$ and $T$ are incomparable if they have different sizes. It follows that the induction of the trivial $kG$-module, corresponding to $\{1,\dots,n\}$, is injective.   
\end{proof}

The asymptotic for the Catalan monoid thus takes the form of \autoref{diagram bn}, (a).

\subsubsection{Full linear monoids}
If $F=\mathbb{F}_{q^s}$ is a finite field of characteristic $q$, we write $\M_m(F)$ for the full linear monoid $\M(m,q^s)$. 
\begin{Proposition}\label{linear monoid inj}
Let $F$ be a field of different characteristic to $k$, then $k\M_m(F)$ satisfies the group-injective condition.
\end{Proposition}
\begin{proof}
If $F$ and $k$ have different characteristics, we know that (see \cite{kovacs1992semigroup}) $k\M_m(F)\cong \prod_{r=0}^m \M_{k_r}(k\GL_r(F))$ where $k_r$ is the number of $r$-dimensional spaces of $F^m$. In particular, $kGL_m(F)$ is a summand of $kM_n(F)$ and so the inductions of projective $kG$-modules are projective $kM$-modules. By \cite[Corollary 15.7]{steinberg2016representation}, if $k$ and $F$ have different characteristics, then $kM_n(F)$ is self-injective. Thus the group-injective condition is satisfied in this case. 
\end{proof}
\begin{Remark}
\autoref{linear monoid inj} settles the `mixed characteristics' case. In the `defining characteristic' case ($\text{char \ } k=\text{char} \ F$), it is known that when $m=2$ the group-injective condition is also satisfied, see \cite[Proposition 4B.3]{he2025growthproblemsrepresentationsfinite}). Note also that the isomorphism $k\M_m(F)\cong \prod_{r=0}^m \M_{k_r}(k\GL_r(F))$ implies that $k\M_m(F)$ is semisimple when $\textup{char} \ k \nmid \lvert \GL_m(F)\rvert$. 
\end{Remark}

This implies an asymptotic formula for the full linear monoid $\M(m,q^s)$ over any field that is not in the defining characteristic. The group of units of $\M(m,q^s)$ is the general linear group, so a closed formula (without characters) is however nasty (unless $m=2$).

\subsection{Temperley--Lieb and Motzkin monoid}

For this section assume $k=\C$, then we have that $\C TL_m$ (resp. $\C Mo_m$) satisfies the group-injective condition if and only if $m\equiv 2 \bmod 3$ (resp. $m\equiv 0 \bmod 2$), see \autoref{P:GI} above. 
\begin{Proposition}\label{TL Mo an} Let $m\ge 1$. If $M=TL_m$ or $Mo_m$ and $V$ is a $\C M$-module with no element apart from 1 acting as identity, then $$a(n)=(\dim V)^n.$$ 
\end{Proposition}
Since the group-injective condition is not always fulfilled, we need a different argument. For a monoid $M$, denote by $1_b$ the induction of the trivial $kG$-module. If $V$ is a $kM$-module, denote by $r(V)$ the number of its composition factors isomorphic to $1_b$ divided by its length. As in \autoref{n0 bound}, let $Z_m$ denote the submodule of $V^{\otimes m}$ on which all nonunits act as 0.
\begin{Proposition} \label{ref finite} \label{composition argument}
Let $M$ be a monoid with trivial group of units. Assume no nonunit act as 1 on $V$. If there is some $c<1$ such that $r(V)<c$ for all indecomposable $kM$-modules $V\not\cong 1_b$, then $b(n)\sim (\dim V)^n$.  
\end{Proposition}
\begin{proof}
Suppose not, then the cell $\Gamma_b^G$ consisting of the single node $1_b$ is never reached. By decomposing $V^{\otimes n}$ into indecomposable summands, we see that $r(V^{\otimes n})\le c$. However, $r(V^{\otimes n})\ge \frac{\dim Z_n}{\dim V^{\otimes n}}$ but the latter tends to 1. To see this, observe that if $v_1\dots v_m$ is a basis of $T$ and all nonunits act as 0 on $v_1$, then all nonunits act as 0 on  $\langle v_{i_1}\otimes v_{i_2}\otimes \dots v_{i_m}\mid i_{j}=1 \text{\ for some\ } j\rangle \le T^{\otimes m}$.   
\end{proof}

\begin{proof}[Proof of \autoref{TL Mo an}]
By \cite[Theorem 3.10]{Bellet_te_2017} both $TL_n$ and $Mo_n$ are representation-finite, and the only indecomposable module with all composition factors isomorphic to $1_b$ is $1_b$ itself. Thus, the hypothesis of \autoref{composition argument} is satisfied. (We note that the dilute Temperley--Lieb algebra $dTL_n(\beta)$ discussed in \cite{Bellet_te_2017} is isomorphic to the Motzkin algebra $Mo_n(\beta +1)$, so the monoid case corresponds to when $\beta = 0$.)  
\end{proof}

\begin{figure}[H]
\centering
\includegraphics[width=0.5\linewidth]{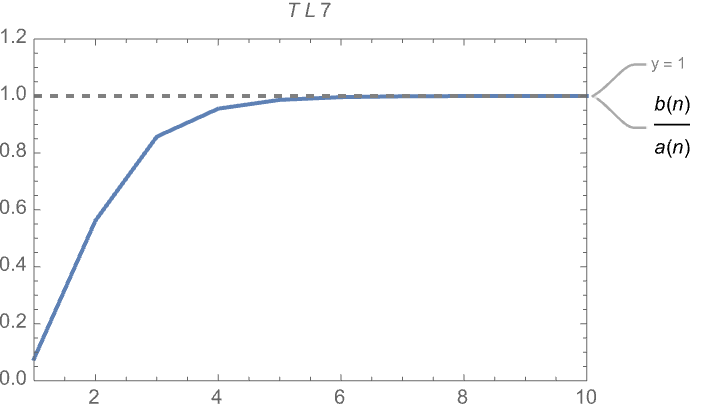}
\caption{The ratio $b(n)/a(n)$ for the 13-dimensional $TL_7$-module $V_3$. In this case $a(n)=13^n.$}
\label{fig:TL7conv}
\end{figure}

We illustrate in \autoref{fig:TL7conv} the ratio $b(n)/a(n)$ when $V=V_3$ and $a(n)=13^n$. In this case the second largest eigenvalue is $4$ and so $\lvert b(n)/a(n)-1\rvert \in \mathcal{O}((4/13)^n)$ and $\lvert b(n)-a(n)\rvert\in \mathcal{O}(4^n)$.

\section{The semigroup case}\label{semigroup section}
If $S$ is not a monoid but a finite semigroup without identity, we can still consider the growth problem associated to a (now non-unital) $kS$-module $V$. Unlike in the monoid case, the \textit{null representation} $Z$ with all elements in $S$ acting as 0 is admitted as a $kS$-module, and it forms a tensor ideal (in the sense that $V\otimes Z=Z^{\oplus \dim V}$ for all $V$.)
We thus obtain the following:

\begin{Proposition}\label{semigroup length}
Let $V$ be a $kS$-module on which no element acts invertibly, then 
$$l(n)\sim k(n)=(\dim V)^n.$$
Moreover, if $n_0$ denotes the smallest natural number such that $[V^{\otimes n_0}:Z]>0$, then $n_0\le L$ where $L$ is the number of $\mathcal{L}$-classes of $S$.
\end{Proposition}
\begin{proof}
Let $Z_n$ denote the submodule of $V^{\otimes n}$ annihilated by all elements in $S$, then by the same argument as in the proof of \autoref{n0 bound}, $Z_L\le V^{\otimes L}$ is nonzero. This implies that the trivial module $k$ (and hence any module, since $V\otimes k\cong V$ and $Z$ is a tensor ideal) has a path to $Z$ in the fusion graph. In other words, $Z$ now plays the same role as $\Gamma^G_l$ does in the monoid setting. This now implies $k(n)=(\dim V)^n$. ($Z$ is a \textit{FBC} in the sense of \cite{lacabanne2024asymptotics}, and by \cite[Proposition 4.22]{lacabanne2024asymptotics} it controls $k(n)$.) 
\end{proof}

We also have the analog of \autoref{conj}:

\begin{Conjecture}\label{semigroup conj}
Let $V$ be a $kS$-module on which no element acts invertibly, then there is some $m\ge 1$ such that $Z$ is a summand of $V^{\otimes m}$, and consequently $$a(n)=(\dim V)^n. $$
\end{Conjecture}
Note that if $Z$ is a summand of $V^{\otimes m}$ for some $m$, then $Z$ plays the same role as $\Gamma_b^G$ in the monoid setting. 

\begin{figure}[H]
\centering
\begin{minipage}{0.2\textwidth}
\centering
\includegraphics[width=\textwidth]{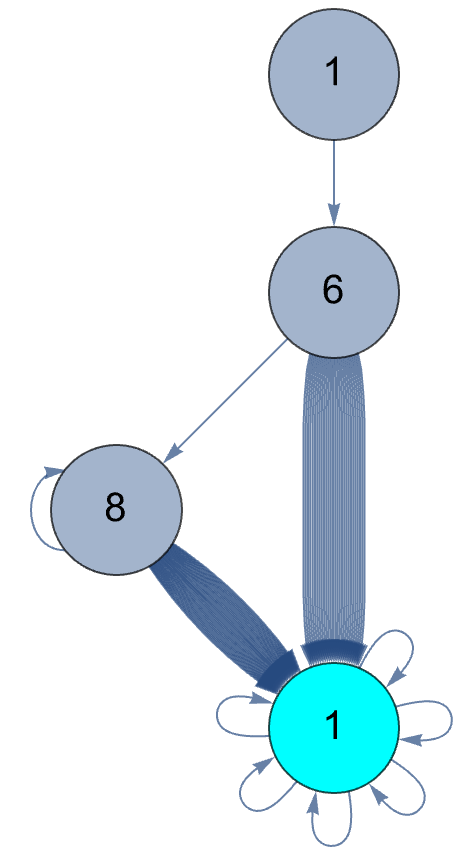} 
\end{minipage}
\begin{minipage}{0.3\textwidth}
\centering     
\[\begin{pmatrix}
0 & 0 & 0 & 0 \\
1 & 0 & 0 & 0 \\
0 & 28 & 6 & 40 \\
0 & 1 & 0 & 1
\end{pmatrix}\]
\end{minipage}
\caption{Here $S$ is the semigroup of nonunits in $Pa_3$, and $V$ is the 6-dimensional indecomposable summand of $kS$. The nodes are labeled with their dimensions and $Z$ is colored in cyan.}
\label{singular ideal}
\end{figure}

\begin{Example}[Singular ideal of partition monoid]
Let $S$ be the semigroup of nonunits in the partition monoid $Pa_3$, and let $V$ be the six-dimensional indecomposable summand of $kS$. In this case $V^{\otimes 2}$ contains $Z$ as a summand, and we have $a(n)=6^n$. We plot the fusion graph and its adjacency matrix in \autoref{singular ideal}.  
\end{Example}

Finally, we give an analog of \autoref{ref finite} in the semigroup setting. For a $kS$-module $V$, let $r(V)$ be the number of its composition factors  isomorphic to $Z$ divided by its length.

\begin{Proposition}
Suppose no element acts invertibly on $V$. If there is $c<1$ such that $r(V)<c$ for all indecomposable $kS$-modules $V \not\cong Z$, then $b(n)=(\dim V)^n.$ 
\end{Proposition}
\begin{proof}
Similar to \autoref{ref finite}.
\end{proof}


\end{document}